\documentclass[12pt]{iopart}

\pdfoutput=1
\usepackage{iopams}
\usepackage{setstack}
\usepackage{hyperref} 
\usepackage{graphicx}
\usepackage{algorithm}

\newtheorem{theorem}{Theorem}
\newcommand{\bb}{\boldsymbol{b}}
\newcommand{\rb}{\boldsymbol{r}}
\newcommand{\ub}{\boldsymbol{u}}
\newcommand{\vb}{\boldsymbol{v}}
\newcommand{\xb}{\boldsymbol{x}}
\newcommand{\wb}{\boldsymbol{w}}
\newcommand{\Ab}{{\mathbf A}}
\newcommand{\Cb}{{\mathbf C}}
\newcommand{\Db}{{\mathbf D}}
\newcommand{\Ib}{{\mathbf I}}
\newcommand{\Lb}{{\mathbf L}}
\newcommand{\Pb}{{\mathbf P}}
\newcommand{\Rb}{{\mathbf R}}

\newcommand{\R}{\mathbb{R}}
\newcommand{\epsb}{\boldsymbol{\epsilon}}
\newcommand{\mub}{\boldsymbol{\mu}}
\newcommand{\thetab}{\boldsymbol{\theta}}

\begin{document}

\title[Semivariogram methods for inverse problems]{Semivariogram methods for modeling Whittle-Mat\'ern priors in Bayesian inverse problems}

\author{Richard D Brown$^1$, Johnathan M Bardsley$^1$ and
Tiangang Cui$^2$}
\address{$^1$ Department of Mathematical Sciences,
University of Montana, Missoula, MT 59812, United States}
\address{$^2$ School of Mathematics, Monash University, Melbourne, Australia}
\eads{\mailto{rick.brown@umontana.edu}, \mailto{bardsleyj@mso.umt.edu}, \mailto{Tiangang.Cui@monash.edu}}
\date{}

\begin{abstract}
We present a new technique, based on semivariogram methodology, for obtaining point estimates for use in prior modeling for solving Bayesian inverse problems. This method requires a connection between Gaussian processes with covariance operators defined by the Mat\'ern covariance function and Gaussian processes with precision (inverse-covariance) operators defined by the Green's functions of a class of elliptic stochastic partial differential equations (SPDEs). We present a detailed mathematical description of this connection. We will show that there is an equivalence between these two Gaussian processes when the domain is infinite -- for us, $\mathbb{R}^2$ -- which breaks down when the domain is finite due to the effect of boundary conditions on Green's functions of PDEs. We show how this connection can be re-established using extended domains. We then introduce the semivariogram method for estimating the Mat\'ern covariance hyperparameters, which specify the Gaussian prior needed for stabilizing the inverse problem. Results are extended from the isotropic case to the anisotropic case where the correlation length in one direction is larger than another. Finally, we consider the situation where the correlation length is spatially dependent rather than constant. We implement each method in two-dimensional image inpainting test cases to show that it works on practical examples.
\end{abstract}


\section{Introduction}
Inverse problems are ubiquitous in science and engineering. They are characterized by the estimation of parameters in a mathematical model from measurements and by a high-dimensional parameter space that typically results from discretizing a function defined on a computational domain. For typical inverse problems, the process of estimating model parameters from measurements is ill-posed, which motivates the use of regularization in the deterministic setting and the choice of a prior probability density in the Bayesian setting. In this paper, we consider linear models of the form
\begin{eqnarray}
\label{b=Ax}
\boldsymbol{b}=\mathbf{A}\boldsymbol{x}+\boldsymbol{\epsilon}, \quad \boldsymbol{\epsilon}\sim\mathcal{N}(\boldsymbol{0},\lambda^{-1}\mathbf{I}_M),
\end{eqnarray}
where $\bb\in\R^M$ is the vector of measurements, $\Ab\in\R^{M\times N}$ is the forward model matrix, $\xb\in\R^N$ is the vector of unknown parameters, and $\epsb\sim\mathcal{N}(\boldsymbol{0},\lambda^{-1}\mathbf{I}_M)$ is the observation noise that follows a zero-mean Gaussian distrubution with covariance matrix $\lambda^{-1}\Ib_M$, 
with $\Ib_M$ denoting the $M\times M$ identity. In typical inverse problems, $\Ab\xb$ is the discretization of a continuous forward model $Ax$, where $A$ is a linear operator and $x$ is a function. The components of the vector $\xb$ satisfy $x_i=x(\ub_i)$, where $\ub_i\in\mathbb{R}^d$ is the location of the $i$th element of the numerical grid. The random vector $\bb$ in (\ref{b=Ax}) has conditional probability density function
\begin{eqnarray}
\label{equ:pdf}
p(\boldsymbol{b}|\boldsymbol{x},\lambda)\propto\exp\left(-\frac{\lambda}{2}\Vert\mathbf{A}\boldsymbol{x}-\boldsymbol{b}\Vert^2   \right),
\end{eqnarray}
where $\propto$ denotes proportionality and $\Vert\cdot\Vert$ denotes the $\ell^2$-norm. 
The maximizer of $p(\bb|\xb,\lambda)$ with respect to $\xb$ is known as the maximum likelihood estimator, and we denote it by  $\xb_{\rm ML}$. As stated above, due to ill-posedness, $\xb_{\rm ML}$ is unstable with respect to errors in $\bb$, i.e., small changes in $\bb$ result in large relative changes in $\xb_{\rm ML}$.

There are various methods to stabilize the solution of inverse problems, but they all involve some form of regularization. In this paper, we take the Bayesian approach \cite{KaiSom}, which requires the definition of a prior probability density function on $\xb$. We make the assumption that the prior is Gaussian of the form $\boldsymbol{x}\sim \mathcal{N}\left(\boldsymbol{0},(\delta\mathbf{P})^{-1}\right)$, which has probability density function
\begin{eqnarray}
\label{equ:generalprior}
p(\boldsymbol{x}|\delta)\propto\exp\left(-\frac{\delta}{2}\boldsymbol{x}^T\mathbf{P}\boldsymbol{x}   \right),
\end{eqnarray}
where $\mathbf{P}$ is the precision (inverse-covariance) matrix. 

Now that we have defined the prior (\ref{equ:generalprior}) and the likelihood (\ref{equ:pdf}), using Bayes' law, we multiply them together to obtain the posterior density function
\begin{eqnarray}
p(\boldsymbol{x}|\boldsymbol{b},\lambda,\delta)&\propto p(\boldsymbol{b}|\boldsymbol{x},\lambda)p(\boldsymbol{x}|\delta)\nonumber\\
&\propto\exp\left(-\frac{\lambda}{2}\Vert\mathbf{A}\boldsymbol{x}-\boldsymbol{b}\Vert^2 -\frac{\delta}{2}\boldsymbol{x}^T\mathbf{P}\boldsymbol{x}\right),\label{posterior}
\end{eqnarray}
whose maximizer, $\xb_{\lambda,\delta}$, is known as the maximum a posteriori (MAP) estimator. The MAP estimator can be equivalently expressed as
\begin{eqnarray*}
\boldsymbol{x}_{\lambda,\delta}&=\mbox{arg}\min_{\boldsymbol{x}}\left\{\frac{\lambda}{2}\Vert \mathbf{A}\boldsymbol{x}-\boldsymbol{b}\Vert^2 +\frac{\delta}{2}\boldsymbol{x}^T\mathbf{P}\boldsymbol{x}  \right\}.
\end{eqnarray*}
Our primary focus in this paper is to provide formulations and hyperparameter selection techniques for prior precision matrices that have an intuitive interpretation and can be used to solve a wide variety of problems.

\subsection{The Mat\'ern Class of Covariance Matrices and Whittle-Mat\'ern Priors}

It remains to define the prior covariance matrix $\mathbf{C}=\Pb^{-1}$. The Mat\'ern class of covariance matrices has garnered much praise \cite{stein2012interpolation} for its flexibility in capturing many covariance structures and its allowance of direct control of the degree of correlation in the vector $\boldsymbol{x}$ \cite{guttorp2006studies}. The Mat\'ern covariance matrix is defined by the Mat\'ern covariance function, which was first formulated by Mat\'ern in 1947 \cite{matern2013spatial},
\begin{equation}
\label{equ:matern}
C(r)=\sigma^2 \frac{(r/\ell)^{\nu}K_{\nu}(r/\ell)}{2^{\nu-1}\Gamma(\nu)},
\end{equation}
where $r$ is the separation distance; $K_\nu(\cdot)$ is the modified Bessel function of the second kind of order $\nu$ \cite{andrews1992special}; $\Gamma(\cdot)$ is the gamma function; $\ell>0$ is the range parameter; $\nu>0$ is the smoothness parameter; and $\sigma^2$ is the marginal variance. Omitting $\sigma^2$ gives the Mat\'ern correlation function. In the isotropic case, when the covariance depends only on the distance between elements, given the covariance parameters $\sigma^2,\nu,$ and $\ell$, one can obtain the covariance matrix $\mathbf{C}$ of a vector $\boldsymbol{x}=[ x_1,\dots,x_N]^T$ with spatial positions $\{ \boldsymbol{u}^T_1, \dots, \boldsymbol{u}^T_N\}\subset\R^d$ by letting
\begin{equation*}
[\mathbf{C}]_{ij}=\mbox{Cov}(x_i,x_j)=C(\Vert \boldsymbol{u}_i-\boldsymbol{u}_j\Vert),
\end{equation*}
where $C$ is defined by (\ref{equ:matern}). 

\begin{figure}
\centering
\includegraphics[width=2.5in]{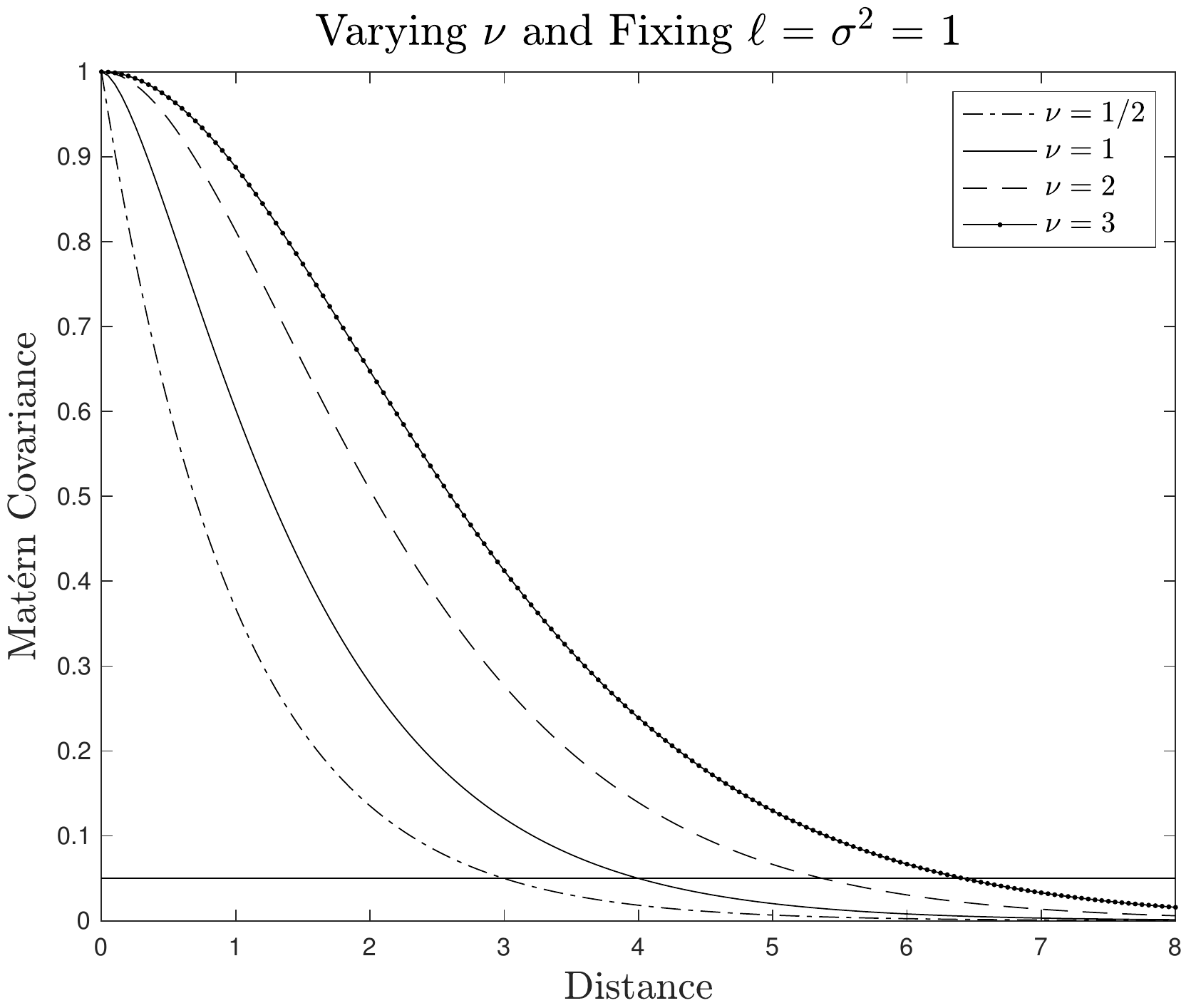}\quad
\includegraphics[width=2.5in]{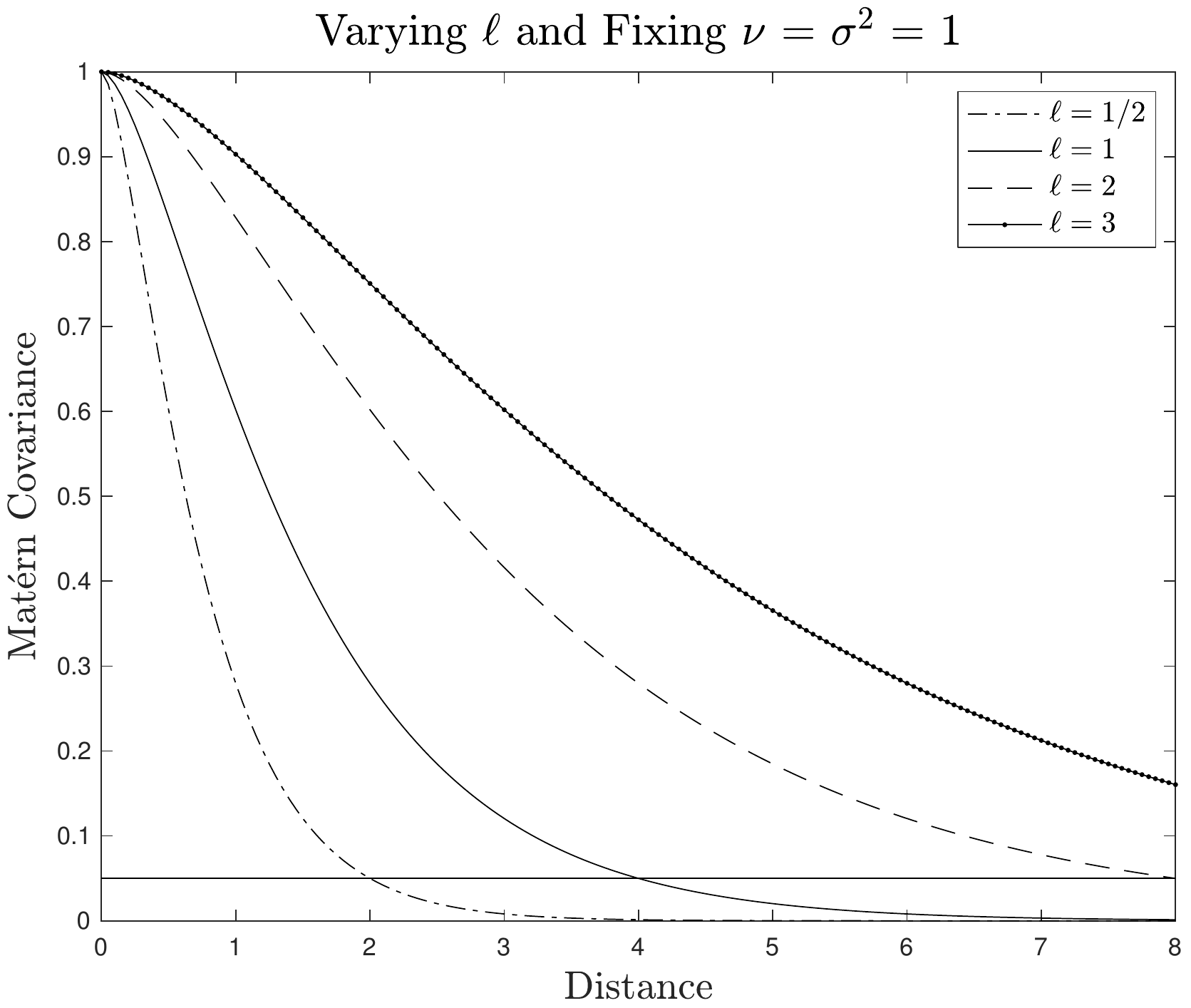}
\caption{Behavior of the Mat\'ern covariance function. The smoothness parameter, $\nu$, primarily affects the covariance at small distances whereas the range parameter, $\ell$, mainly affects the decay rate of the covariance. The horizontal line corresponds to a covariance value of 0.05 and the practical range is the distance at which the covariance intersects this line.}
\label{materncovs}
\end{figure}

The parameters of the Mat\'ern covariance function are not as straightforward to interpret as the parameters of some other covariance functions. When $\nu$ is small ($\nu\rightarrow 0^+$), the spatial process is said to be rough, and when it is large ($\nu\rightarrow \infty$), the process is smooth \cite{guttorp2006studies,minasny2005matern}. Figure \ref{materncovs} shows how the covariance function behaves with different values of $\ell$ and $\nu$: on the left, $\ell=\sigma^2=1$ and $\nu$ varies, while on the right $\nu=\sigma^2=1$ and $\ell$ varies. Note that as $\nu$ increases, the behavior at small lags changes, leading to more correlation at smaller distances and a larger {\em practical range}, which is defined to be the distance at which the correlation is equal to 0.05. In Figure 1, this is the distance at which the covariance function intersects the horizontal line.  Meanwhile, as $\ell$ decreases, the decay rate of the covariance increases considerably, which decreases the practical range. Although $\ell$ is known as the range parameter, the parameter $\nu$ also affects the practical range. In \cite{lindgren2011explicit}, a range approximation $\rho=\ell\sqrt{8\nu}$ is used where $C(\rho)\approx 0.10$.

Despite the benefits of using the Mat\'ern class of covariance matrices, its use can be problematic for inverse problems because computing the precision matrix $\Pb$, which is what appears in the posterior (\ref{posterior}), requires inverting a dense $N\times N$ matrix. Using the fast Fourier transform (FFT) \cite{wood1994simulation,dietrich1997fast,bardsley2018} to operate with $\Pb$ and $\Cb$ more efficiently is recommended if $\xb$ is defined on a regular grid and periodic boundary conditions are assumed. 
In other cases, it is useful that the Mat\'ern covariance function has a direct connection to a class of elliptic SPDEs \cite{lindgren2011explicit} whose numerical discretization yields sparse precision matrices, $\Pb$, that are computationally feasible to work with even when $N$ is large. Connections of this type were first shown to exist by Whittle in \cite{whittle1954stationary}, where he showed the connection held for a special case of the Mat\'ern covariance class. Hence, priors that depend on this connection are often referred to as Whittle-Mat\'ern priors. The connection between the general Mat\'ern covariance function and SPDEs has been used in a wide range of applications for defining computationally feasible priors for high-dimensional problems \cite{roininen2016hyperpriors, roininen2014whittle, monterrubio2018posterior}. Moreover, work has been done in establishing convergence theorems for, and lattice approximations of, these Whittle-Mat\'ern priors \cite{roininen2013constructing}.

The remainder of the paper is organized as follows. In Section \ref{sec:Matern}, we describe in detail the connection between zero-mean Gaussian processes with the isotropic Mat\'ern covariance operator and those that arise as solutions of a class of elliptic SPDEs. In Section \ref{sec:variogram}, we show how to estimate the hyperparameters in the isotropic Whittle-Mat\'ern prior using the semivariogram method, and then we show how to use this approach to define the prior when solving a Bayesian inverse problem. In Section \ref{sec:anis}, we extend these ideas to the anisotropic case and then we consider images with regions that require different covariance structures in Section \ref{sec:regional}. For each section, we present numerical tests on two-dimensional image inpainting test cases. We end with conclusions in Section \ref{sec:conclusions}.

\section{Whittle-Mat\'ern Class Priors via SPDEs}
\label{sec:Matern}

In this section, we will show that the Whittle-Mat\'ern class of priors can be specified as the solution of the SPDE
\begin{equation}
\label{equ:SPDE}
(1-\ell^{2}\Delta )^{\beta/2}x(\boldsymbol{u}) =\mathcal{W}(\boldsymbol{u}), \quad \boldsymbol{u}\in\mathbb{R}^d,\quad \beta=\nu+d/2, \quad \ell,\nu>0,
\end{equation}
where $\Delta=\sum_{i=1}^d\frac{\partial^2}{du_i^2}$ is the Laplacian operator in $d$ dimensions, and $\mathcal{W}$ is spatial Gaussian white noise with unit variance, which we define below. Although this connection has been shown to exist \cite{whittle1954stationary, lindgren2011explicit, roininen2014whittle}, here we provide a significantly more detailed derivation of this result than we have seen elsewhere. Our derivation is based on the Green's function of the differential operator. For other linear operators with sufficient smoothness, e.g., the one in the Stokes equations and the one in the heat equation, the corresponding SPDEs can be used to define different Gaussian processes \cite{lord2014introduction}. The method we employ here provides a potential way to derive the covariance functions of the Gaussian processes induced by other linear SPDEs as well.

\subsection{Preliminary Definitions}
Before deriving the solution of (\ref{equ:SPDE}), we need some preliminary definitions.

\subsubsection{Gaussian Fields}

A stochastic process $\{x(\boldsymbol{u}), \boldsymbol{u}\in \Omega\}$, with $\Omega\subset\mathbb{R}^d$, is a {\em Gaussian field}  \cite{rue2005gaussian} if for any $k\ge1$ and any locations $\boldsymbol{u}_1,\dots,\boldsymbol{u}_k\in \Omega$, $[x(\boldsymbol{u}_1),\dots,x(\boldsymbol{u}_k)]^T$ is a normally distributed random vector with mean $\mub=\big[E[x(\boldsymbol{u}_1)],\ldots,E[x(\ub_k)]\big]^T$, where $E[\hspace{.5mm}\cdot\hspace{.5mm}]$ denotes expected value, and covariance matrix $[\Cb]_{ij}= \mbox{Cov}(x(\boldsymbol{u}_i),x(\boldsymbol{u}_j))=E[(x(\boldsymbol{u}_i)-E[x(\ub_i)])(x(\boldsymbol{u}_j)-E[x(\ub_j)])]$, for $1\leq i,j\leq k$. The covariance function is defined $C(\boldsymbol{u}_i,\boldsymbol{u}_j):= \mbox{Cov}(x(\boldsymbol{u}_i),x(\boldsymbol{u}_j))$. It is necessary that the covariance function is positive definite, i.e., for any $\{\boldsymbol{u}_1,\dots,\boldsymbol{u}_k\}$, with $k\ge1$, the covariance matrix $\Cb$ defined above
is positive definite. The Gaussian field is called {\em stationary} if the mean is constant and the covariance function satisfies $C(\ub,\vb)=C(\ub-\vb)$ and {\em isotropic} if $C(\boldsymbol{u},\boldsymbol{v})=C(\Vert\ub-\vb\Vert)$.

\subsubsection{White Noise}

The term {\em white noise} \cite{lord2014introduction,walsh1986introduction} comes from light. White light is a homogeneous mix of wavelengths, as opposed to colored light, which is a heterogeneous mix of wavelengths. In a similar way, white noise contains a homogeneous mix of all the different basis functions. The mixing of these basis functions is determined by a random process. When this random process is Gaussian, we have {\em Gaussian white noise}. Consider a domain $\Omega$ and let $\{\phi_j:j=1,2,\dots\}$ be an orthonormal basis of $L^2(\Omega)$ where $L^2(\Omega)=\left\{f:\Omega\rightarrow \mathbb{R}\mid\int_\Omega|f(x)|^2dx<\infty\right\}$. Then Gaussian white noise is defined by
\begin{equation}
\label{equ:white}
\mathcal{W}(\boldsymbol{u})=\sum_{j=1}^\infty\xi_j\phi_j(\boldsymbol{u}),\quad \xi_j\stackrel{iid}{\sim}\mathcal{N}(0,\eta^2).
\end{equation}
If we are dealing with spatial Gaussian white noise with unit variance, then $\boldsymbol{u}$ refers to location and $\eta^2=1$. With this definition, it is clear that Gaussian white noise has mean zero:
$
E[\mathcal{W}(\boldsymbol{u})]=\sum_{j=1}^\infty E\left[\xi_j\right]\phi_j(\boldsymbol{u})=0.
$
Moreover, one can show that
$
\mbox{Cov}\left(\mathcal{W}(\boldsymbol{u}),\mathcal{W}(\boldsymbol{v})\right)=\eta^2\delta_f(\boldsymbol{u}-\boldsymbol{v}),
$
where $\delta_f(\cdot)$ is the {\em Dirac delta function} \cite{hassani2009dirac}, also known as the delta distribution. We include the subscript $f$ to differentiate the delta function from the $\delta$ hyperparameter used elsewhere in this paper. A well-known and very important property of the Dirac delta function is that it satisfies the {\em sifting property}:
$
f(\ub)=\int_{\R^d}\delta_f(\ub-\vb)f(\vb)d\vb.
$

\subsubsection{Green's Functions}
\label{sec:Green}

We now consider differential equations of the form $Lx(\ub)=f(\ub)$, $\ub\in\R^d$, where $L$ is a linear, differential operator. A {\em Green's function} \cite{gockenbach2005partial, stakgold2011green}, $g$, of $L$ is any solution of
$
Lg(\ub,\vb)=\delta_f(\ub-\vb)
$.
Using the Green's function, the solution of the equation $Lx(\ub)=f(\ub)$ can be written as
\begin{equation}
\label{equ:green}
x(\ub)=\int_{\R^d} g(\ub,\vb)f(\vb)d\vb.
\end{equation}

\subsection{The Gaussian Field Solution of the SPDE (\ref{equ:SPDE})}
\label{sec:iso_cov}

In this subsection, we will prove the following theorem concerning the solution of the SPDE (\ref{equ:SPDE}).

\begin{theorem}
\label{thm:connection}
The solution $x(\boldsymbol{u})$ of (\ref{equ:SPDE}) is a Gaussian field with mean zero and Mat\'ern covariance function defined by (\ref{equ:matern}).
\end{theorem}

\noindent {\em Proof.} To begin, we note that the Green's function for (\ref{equ:SPDE}) is the solution of
\begin{equation}
\label{eq:PDE=dirac}
(1-\ell^{2}\Delta )^{\beta/2}g(\boldsymbol{u},\boldsymbol{v}) = \delta_f(\boldsymbol{v}-\boldsymbol{u}).
\end{equation}
Using (\ref{equ:green}), the solution to (\ref{equ:SPDE}) is given by
\begin{equation}
\label{equ:SPDEsolution}
x(\boldsymbol{u})=\int_{\mathbb{R}^d} g(\boldsymbol{u},\boldsymbol{v})\mathcal{W}(\boldsymbol{v})d\boldsymbol{v},
\end{equation}
making $x(\boldsymbol{u})$ a Gaussian field since it is a linear transformation of Gaussian white noise. 

We now compute the mean and covariance of the Gaussian field, $x(\ub)$, defined by (\ref{equ:SPDEsolution}). Since the Green's function is a strictly-positive, symmetric, and rapidly decaying function, we can apply Fubini's theorem \cite{saks1937theory} to obtain the mean of $x(\ub)$:
$$
E[x(\boldsymbol{u})]=E\left[ \int_{\mathbb{R}^d} g(\boldsymbol{u},\boldsymbol{v})\mathcal{W}(\boldsymbol{v})d\boldsymbol{v}  \right]
= \int_{\mathbb{R}^d} g(\boldsymbol{u},\boldsymbol{v})E\left[\mathcal{W}(\boldsymbol{v}) \right]d\boldsymbol{v}
=0.
$$
Since $x(\ub)$ has mean zero, the covariance is given by
\begin{eqnarray*}
\mbox{Cov}(x(\boldsymbol{u}),x(\boldsymbol{u'}))&=E[x(\boldsymbol{u})x(\boldsymbol{u'})]\\
&=\int_{\mathbb{R}^d}\left(\int_{\mathbb{R}^d} E[\mathcal{W}(\boldsymbol{v})\mathcal{W}(\boldsymbol{v'})]g(\boldsymbol{u},\boldsymbol{v})d\boldsymbol{v}\right) g(\boldsymbol{u'},\boldsymbol{v'})d\boldsymbol{v'}\\
&=\int_{\mathbb{R}^d}\left(\int_{\mathbb{R}^d} \delta_f(\boldsymbol{v}-\boldsymbol{v'})g(\boldsymbol{u},\boldsymbol{v})d\boldsymbol{v} \right)g(\boldsymbol{u'},\boldsymbol{v'})d\boldsymbol{v'}\\
&=\int_{\mathbb{R}^d} g(\boldsymbol{u},\boldsymbol{v'}) g(\boldsymbol{u'},\boldsymbol{v'})d\boldsymbol{v'}.
\end{eqnarray*}
If we define $C(\boldsymbol{u},\boldsymbol{u'}):=\mbox{Cov}(x(\boldsymbol{u}),x(\boldsymbol{u'}))$, the previous result implies that if $L=(1-\ell^{2}\Delta)^{\beta/2}$, then for our linear $L$ acting only on $\boldsymbol{u'}$,
\begin{eqnarray}
LC(\boldsymbol{u},\boldsymbol{u'})&=L\int_{\mathbb{R}^d} g(\boldsymbol{u},\boldsymbol{v'}) g(\boldsymbol{u'},\boldsymbol{v'})d\boldsymbol{v'}\nonumber\\
&=\int_{\mathbb{R}^d} \bigg[Lg(\boldsymbol{u'},\boldsymbol{v'})\bigg]g(\boldsymbol{u},\boldsymbol{v'})d\boldsymbol{v'}\nonumber \\
&=\int_{\mathbb{R}^d} \delta_f(\boldsymbol{u'}-\boldsymbol{v'}) g(\boldsymbol{u},\boldsymbol{v'})d\boldsymbol{v'}\nonumber\\
&=g(\boldsymbol{u},\boldsymbol{u'}).\label{Lcov}
\end{eqnarray}

To derive the Green's function $g$ in (\ref{Lcov}), we first define $g(\boldsymbol{u}):=g(\boldsymbol{u},\boldsymbol{0})$. Then (\ref{eq:PDE=dirac}) implies
\begin{equation}
\label{equ:SPDEgreen}
(1-\ell^{2}\Delta )^{\beta/2}g(\boldsymbol{u}) = \delta_f(\boldsymbol{u}).
\end{equation}
To proceed, we must take the Fourier transform \cite{sneddon1995fourier,kwasnicki2017ten} of both sides of (\ref{equ:SPDEgreen}). This yields
\begin{eqnarray*}
(1+\ell^{2}\Vert\boldsymbol{\omega}\Vert^2 )^{\beta/2}\hat{g}(\boldsymbol{\omega}) &=1,
\end{eqnarray*}
where $\boldsymbol{\omega}\in \mathbb{C}^d$ are the coordinates in the Fourier-transformed space and the hat ($\hat{f}$) notation denotes the Fourier-transform of a function $f$. Thus, the Fourier transform of the Green's function is
\begin{equation}
\label{equ:greenSPDE}
\hat{g}(\boldsymbol{\omega})=(1+\ell^{2}\Vert\boldsymbol{\omega}\Vert^2)^{-\beta/2}.
\end{equation}

Next, we assume {\em stationarity} so that the covariance only depends on the relative locations of the points, i.e., $\boldsymbol{r}:= \boldsymbol{u}-\boldsymbol{v}$. Then $E[x(\boldsymbol{u})x(\boldsymbol{v})]=E[x(\boldsymbol{r})x(\boldsymbol{0})]=C(\boldsymbol{r},\boldsymbol{0}):= C(\boldsymbol{r})$ and (\ref{Lcov}) can be expressed
$
LC(\boldsymbol{r})=g(\boldsymbol{r}).
$
If we take the Fourier transform of both sides of this equation, and appeal to (\ref{equ:greenSPDE}), we obtain
\begin{eqnarray*}
\hat{C}(\boldsymbol{\omega})=(1+\ell^{2}\Vert\boldsymbol{\omega}\Vert^2)^{-\beta}.
\end{eqnarray*}
Since the Laplacian, $\Delta$, is invariant under rotations and translations, we have radial symmetry, which is analogous to isotropy in the covariance. Thus we can let $s=\Vert\boldsymbol{\omega}\Vert$ and $r=\Vert\boldsymbol{r}\Vert$ to obtain the equivalent expression
\begin{equation}
\label{fourier_cov}
\hat{C}(s)=(1+\ell^{2}s^2)^{-\beta}.
\end{equation}
To transform back to the original ($r$) space, we use the Hankel transform \cite{piessens2000hankel} and its relationship to the radially symmetric Fourier transform, i.e.,
\begin{equation*}
\label{HankelFourier}
s^{\frac{d-2}{2}}\hat{C}(s)=(2\pi)^{\frac{d}{2}}\int_{0}^\infty J_{\frac{d-2}{2}}(sr)r^{\frac{d-2}{2}} C(r)rdr,
\end{equation*}
where $C$ is the original (untransformed) covariance function and $J_\nu(\cdot)$ is the Bessel function of the first kind of order $\nu$; see \cite[Section 2]{grafakos2013fourier} for a proof.
Using appropriate substitutions in the inverse Hankel transform and (\ref{fourier_cov}), we obtain
\begin{eqnarray*}
C(r) =\frac{(2\pi)^{-\frac{d}{2}}}{r^{\frac{d-1}{2}} }\int_0^\infty  J_{\frac{d-2}{2}}(sr) s^{\frac{d-1}{2}}(1+\ell^{2}s^2)^{-\beta}(sr)^{1/2}ds.
\end{eqnarray*}
Finally, using the integral identity \cite[Eq. 20, p. 24, vol. II]{bateman1954tables} and some algebra, we obtain
\begin{eqnarray}
C(r)
&=\frac{\ell^{-\beta-\frac{d}{2}}r^{\beta-\frac{d}{2}}K_{\frac{d}{2}-\beta}(r/\ell)}{(2\pi)^{\frac{d}{2}}2^{\beta-1}\Gamma(\beta)}.\label{C temp}
\end{eqnarray}
Using the fact that $K_{\nu}=K_{-\nu}$, and defining $\sigma^2:=\Gamma(\nu)[\ell^d(4\pi)^{d/2}\Gamma(\nu+d/2)]^{-1}$ with $\nu:= \beta-d/2$, it can be shown that (\ref{C temp}) is exactly the Mat\'ern covariance function (\ref{equ:matern}). 

\hfill $\square$

\subsection{The Effect of a Finite Domain and Boundary Conditions}
\label{boundary_effects}

The proof of Theorem \ref{thm:connection} above assumed that the domain was all of $\mathbb{R}^d$, i.e. $\Omega=\mathbb{R}^d$. However, when solving inverse problems, $x(\boldsymbol{u})$ is restricted to a finite domain $\Omega\subset \mathbb{R}^d$. In such cases, boundary conditions that modify the Green's function must be assumed, and thus the equivalence between the Gaussian fields defined by the SPDE (\ref{equ:SPDE}) and those defined by the Mat\'ern covariance function may not hold.

To see this, consider the case where $d=2$ and $\Omega=[0,1]\times[0,1]$ with Dirichlet (zero) boundary conditions, $x(0,t)=x(1,t)=x(s,0)=x(s,1)=0$, where $0\leq s,t\leq 1$. Additionally, we assume $\nu=1$ so that the exponent of the differential operator is equal to one, making the discretization straightforward. In this case,  (\ref{equ:SPDE}) simplifies to
\begin{equation*}
(1-\ell^2\Delta)x(\ub)=\mathcal{W}(\ub),\quad \ub\in\mathbb{R}^2,\quad \ell>0.
\end{equation*}
Using a uniform mesh on $[0,1]\times[0,1]$ with a step size of $h=1/n$, so that $N=n^2$, yields the numerical discretization
\begin{equation*}
(\mathbf{I}_N+(\ell/h)^2\mathbf{L}_{\mbox{\scriptsize2D}})\boldsymbol{x}=\delta^{-1/2}\boldsymbol{\xi},\quad \boldsymbol{\xi}\sim \mathcal{N}(\boldsymbol{0},\boldsymbol{I}_n),
\end{equation*}
where $\delta$ is the scaling parameter for the prior and $(1/h^2)\Lb_{\mbox{\scriptsize2D}}\xb$ is the standard finite-difference discretization of 
$\left(-\partial^2x(\ub)/\partial u_1^2-\partial^2x(\ub)/\partial u_2^2\right)$ \cite{bardsley2018}.
Then the probability density for $\xb$ is given by
\begin{equation*}
\boldsymbol{x}|\delta,\ell\sim \mathcal{N}\left(\boldsymbol{0},\delta^{-1}(\mathbf{I}+(\ell/h)^2\mathbf{L}_{\mbox{\scriptsize2D}})^{-2}\right),
\end{equation*}
or equivalently,
\begin{equation}
\label{equ:prior_nu1}
p(\boldsymbol{x}|\delta,\ell)\propto\exp\left(-\frac{\delta}{2}\boldsymbol{x}^T(\mathbf{I}+(\ell/h)^2\mathbf{L}_{\mbox{\scriptsize2D}})^2\boldsymbol{x}\right).
\end{equation}
When discretizing the SPDE, there is a scaling factor needed that guarantees that the variance scales systematically with respect to the change of the length-scaling parameter, $\ell$. The exact form of this scaling factor is unimportant for our purposes since we are ultimately only interested in a regularization parameter, $\alpha$, as will be seen in Section \ref{sec:MAP}. To keep notation simpler, we use $\delta$ as a placeholder for this term. This is also the reason we are interested in whether the Mat\'ern correlation rather than the covariance is preserved when restricting our Gaussian field to a finite domain. 

\begin{figure}
\centering
\includegraphics[width=1.7in]{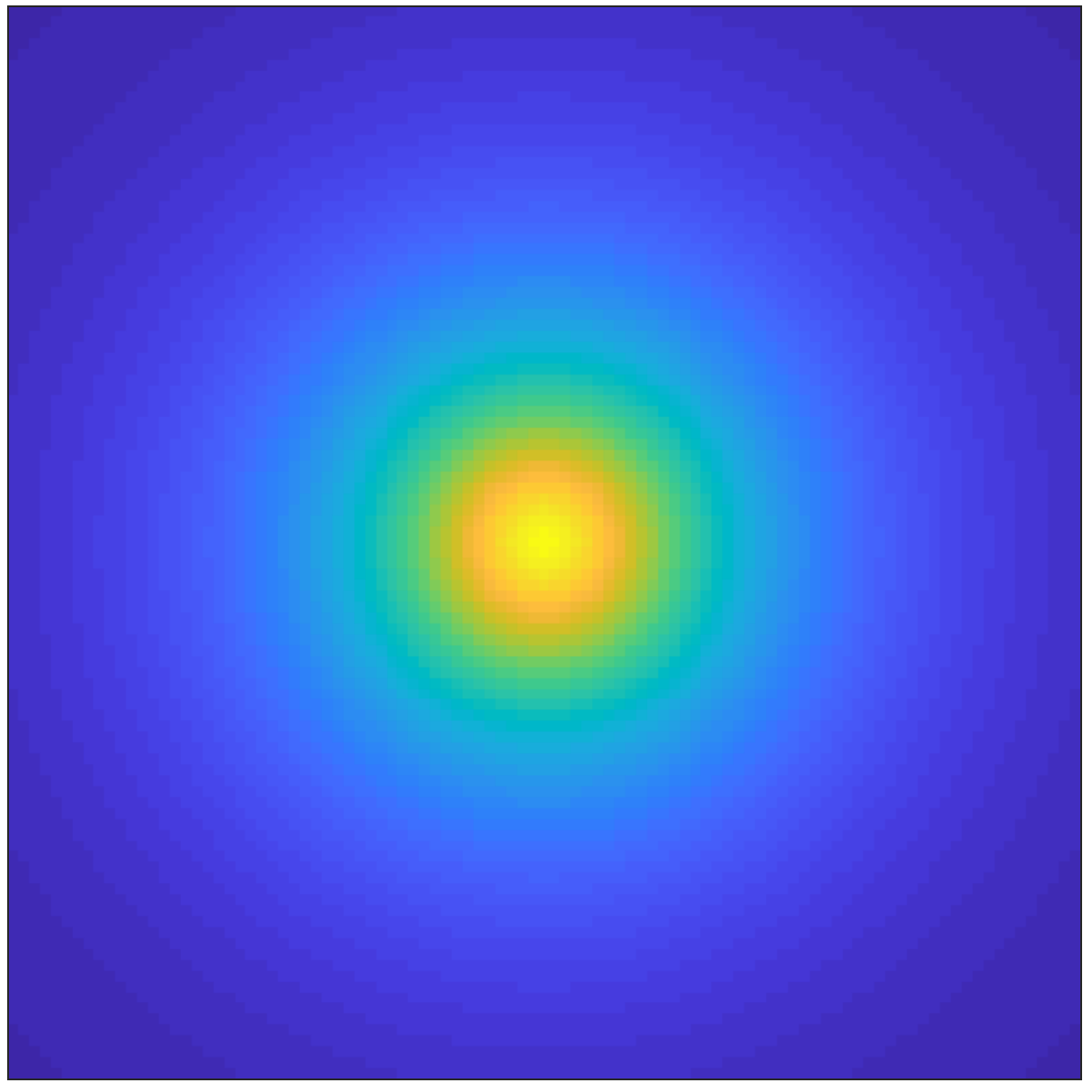}
\includegraphics[width=1.7in]{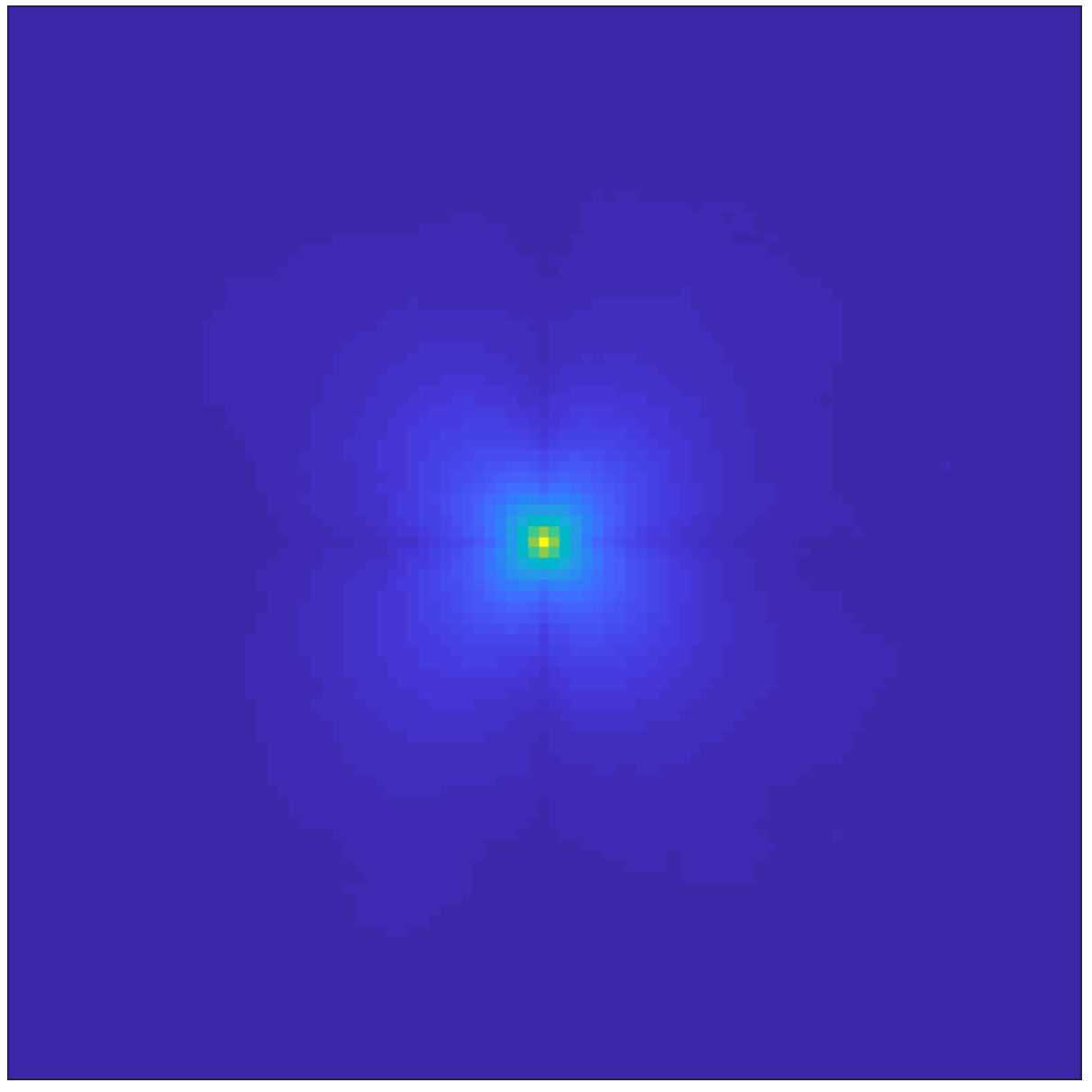}
\includegraphics[width=1.7in]{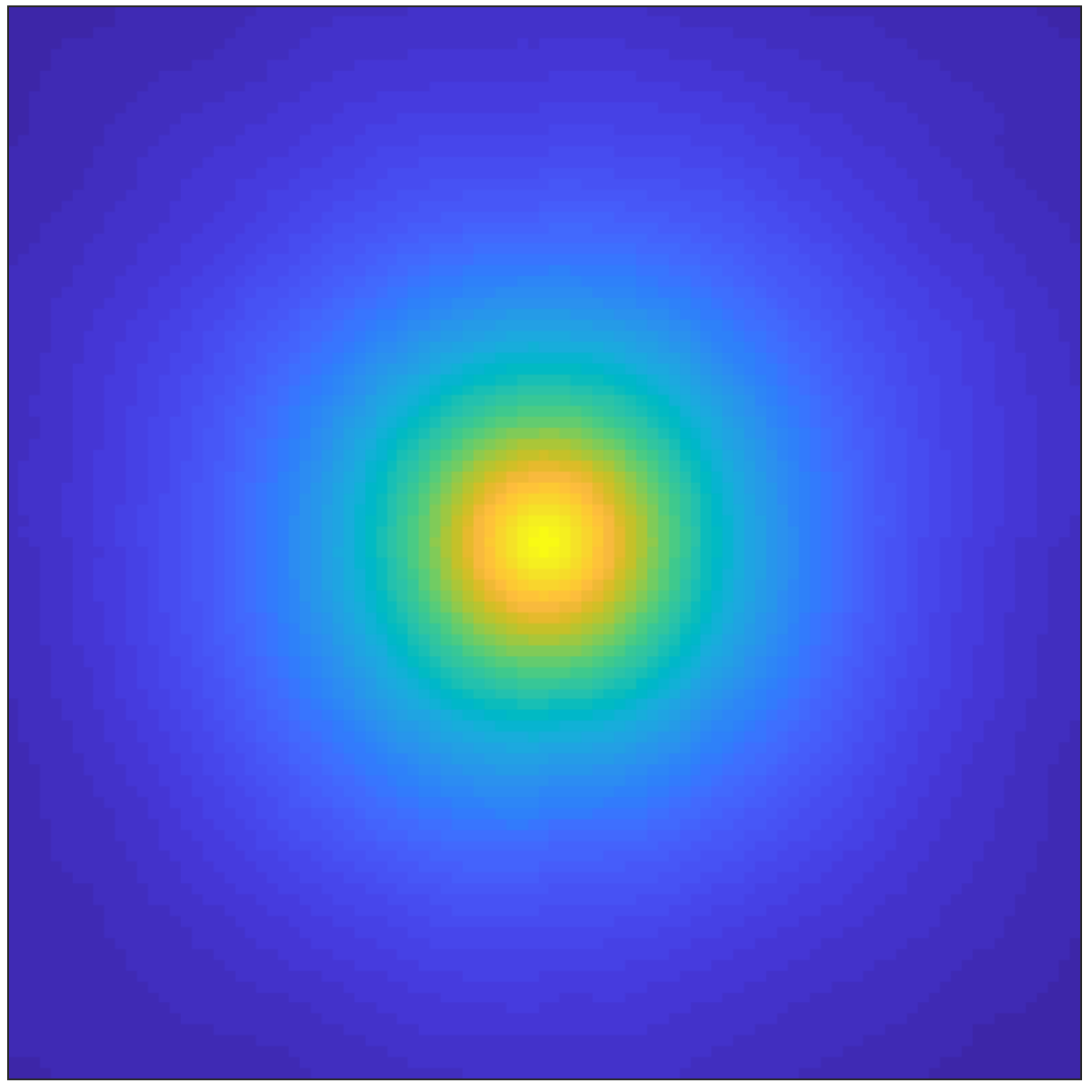}
\caption{Isotropic correlation maps. Plots of the Mat\'ern correlation map (left), the empirical correlation map with $n=50$ computed on the domain $\Omega=[0,1]\times[0,1]$ (middle), and the empirical correlation map computed on the domain $\overline{\Omega}=[-0.5,1.5]\times[-0.5,1.5]$ (right), computed from random draws from the prior (\ref{equ:prior_nu1}) in 2D with $\nu = 1$ and $\ell = 1/4$.}
\label{iso_disconnect}
\end{figure}

We now let $n=50$, so $N=50^2=2500$, and generate $50\,000$ samples from (\ref{equ:prior_nu1}) for each of $N$ $x_{i}$ values, calculate the empirical correlation between the samples, and compare this with the theoretical correlation defined by the Mat\'ern covariance function. We do this for $\ell=1/4$ and plot the results in the middle of Figure \ref{iso_disconnect}, together with the Mat\'ern correlation map on the left. It is clear that there is a disconnection between the empirical correlation and the Mat\'ern correlation.

It is crucial that the connection between the Gaussian fields defined by the SPDE and those defined by the Mat\'ern covariance function holds because then the parameters in the SPDE can be estimated using the semivariogram method described in Section \ref{sec:variogram}. Fortunately, we can restore this connection by extending the computational domain. In two dimensions, we define $\overline\Omega=[1-a,a]\times[1-a,a]$, for $a>1$, e.g., if $a=1.5$ then $\overline\Omega=[-0.5,1.5]\times[-0.5,1.5]$. We then generate realizations for $((2a-1)n)^2=(2n)^2=10\,000$ $x_i$ values on the extended domain
and compute the empirical correlation only for the $x_i$ values that correspond to the original domain, $\Omega=[0,1]\times[0,1]$. The results are plotted on the right side of Figure \ref{iso_disconnect}, where it is clear that the empirical correlation map is nearly indistinguishable from those obtained using the Mat\'ern correlation function.

To determine the $a$ value that extends the domain far enough to restore the Mat\'ern/SPDE connection, but not so far as to introduce unnecessary computational cost, we look to the Mat\'ern correlation function itself. We want to extend the domain far enough so that all $\boldsymbol{x}$ values in $[0,1]\times[0,1]$ have a sufficiently low correlation with the $\boldsymbol{x}$ values at the end of the extended domain. The criterion we used to determine if the connection was restored was based on relative error: $\Vert \boldsymbol{\rho}-\boldsymbol{\rho}_a\Vert_F/\Vert\boldsymbol{\rho}\Vert_F<0.05$, where $\boldsymbol{\rho}$ is the true Mat\'ern correlation matrix, $\boldsymbol{\rho}_a$ is the approximate correlation matrix obtained by discretizing the SPDE, and $\Vert\cdot\Vert_F$ denotes the Frobenius norm.

In tests, it was found that we should always extend the domain at least slightly. If we let $r_{c}$ be the distance for which the Mat\'ern correlation is approximately equal to $c$, then our tests showed that setting $a=1+r_{0.30}$ restores the connection to the Mat\'ern covariance for $\nu\ge1/2$ when using zero boundary conditions and setting $a=1+r_{0.20}$ restores the connection to the Mat\'ern covariance for $\nu\ge1/2$ when periodic boundary conditions are used. For $\nu=1$ and $\ell=1/4$, $a$ should be set to $1.5$ in the Dirichlet boundary condition case, which gives a relative error in the difference of the correlation matrices of 0.0375, and it should be set to $1.6$ when using periodic boundary conditions. We note that since $\ell$ is directly related to the degree of correlation in the prior, the extension necessary to preserve the connection rises sharply as $\ell$ increases. 
It is rare in practice, however, to have $\ell\ge1/4$ when $\nu\ge1$ since that implies the correlation persists across the entire region. Thus, it is uncommon to have to extend beyond a domain of $[-0.5,1.5]\times[-0.5,1.5]$.

For the above discussion, we focused on zero boundary conditions. Similar results hold if periodic boundary conditions are assumed, in which case $\Lb$, and thus $\Lb_{\mbox{\scriptsize2D}}$, can be diagonalized by the FFT, assuming $\xb$ is defined on a regular grid. The FFT-based diagonalization of $\Lb_{\mbox{\scriptsize2D}}$ can be exploited to greatly reduce computational cost, thus when extending the domain in two-dimensions, it is advantageous to use periodic boundary conditions and the extended domain $\overline\Omega=[-0.5,1.5]\times[-0.5,1.5]$ so that $\Lb_{\mbox{\scriptsize2D}}$ defined on $\overline\Omega$ can be diagonalized by the FFT. A more thorough description of the effects of boundary artifacts with different boundary conditions can be found in  \cite{khristenko2019analysis}.

Finally, in our numerical experiment above, we chose a specific value of $\nu$, but other values of $\nu$ can be chosen. The general form of the isotropic prior density in two dimensions, with $\nu$ included as a hyperparameter, is
\begin{equation}
\label{equ:prior}
p(\boldsymbol{x}|\delta,\nu,\ell)\propto\exp\left(-\frac{\delta}{2}\boldsymbol{x}^T(\mathbf{I}+(\ell/h)^2\Lb_{\mbox{\scriptsize2D}})^{\nu+d/2}\boldsymbol{x}\right).
\end{equation}
If $\nu+d/2$ is a non-integer, a fractional power of $\mathbf{I}+(\ell/h)^2\mathbf{L}_{\mbox{\scriptsize2D}}$ must be computed, which is possible, generally speaking, if we have a diagonalization of $\mathbf{I}+(\ell/h)^2\mathbf{L}_{\mbox{\scriptsize2D}}$ in hand, but the resulting precision matrix is typically full and dense. Such a diagonalization is typically computable in one-dimensional examples, even with dense matrices. In two dimensions, however, an efficient diagonalization is possible only if periodic boundary conditions are assumed. We will restrict the exponent $\nu+d/2$ to be an integer in this paper to preserve the sparsity in the precision matrix, which will be especially useful in Section \ref{sec:regional}.

\subsection{Computing MAP Estimators for Whittle-Mat\'ern Priors}
\label{sec:MAP}
Using Bayes' law, we multiply the prior (\ref{equ:prior}) by the likelihood (\ref{equ:pdf}) to obtain the posterior density function
\begin{eqnarray*}
p(\boldsymbol{x}|\boldsymbol{b},\lambda,\delta,\nu,\ell)&\propto& p(\boldsymbol{b}|\boldsymbol{x},\lambda)p(\boldsymbol{x}|\delta,\nu,\ell)\\
&\propto&\exp\left(-\frac{\lambda}{2}\Vert \mathbf{A}\boldsymbol{x}-\boldsymbol{b}\Vert^2-\frac{\delta}{2}\boldsymbol{x}^T(\mathbf{I}+(\ell/h)^2\Lb_{\mbox{\scriptsize2D}})^{\nu+d/2}\boldsymbol{x} \right).
\end{eqnarray*}
The maximizer of $p(\xb|\bb,\lambda,\delta,\nu,\ell)$ is known as the MAP estimator, and it can be computed by solving
\begin{eqnarray}
\boldsymbol{x}_\alpha
&=\mbox{arg}\min_{\boldsymbol{x}}\left\{\frac{1}{2}\Vert \mathbf{A}\boldsymbol{x}-\boldsymbol{b}\Vert^2 +\frac{\alpha}{2}\boldsymbol{x}^T(\mathbf{I}+(\ell/h)^2\mathbf{L}_{\mbox{\scriptsize2D}})^{\nu+d/2}\boldsymbol{x}\right\}\nonumber\\
&=\left(\mathbf{A}^T\mathbf{A}+\alpha(\mathbf{I}+(\ell/h)^2\mathbf{L}_{\mbox{\scriptsize2D}})^{\nu+d/2}\right)^{-1}\mathbf{A}^T\boldsymbol{b},\label{equ:MAP}
\end{eqnarray}
where $\alpha=\delta/\lambda$. Assuming we know $\ell$ and $\nu$, $\alpha$ can be estimated using one of many regularization parameter selection methods (see, e.g.,\cite{vogel2002computational, hansen2005rank,bardsley2018}). One such method is generalized cross validation (GCV): \begin{equation}
\label{equ:GCV}
\alpha=\mbox{arg}\min_{\eta>0}\left\{\frac{\left\Vert \mathbf{A}\Big(\mathbf{A}^T\mathbf{A}+\eta\Pb\Big)^{-1}\mathbf{A}^T\boldsymbol{b}-\boldsymbol{b}\right\Vert^2}{\mbox{tr}\bigg(\mathbf{I}- \mathbf{A}\Big(\mathbf{A}^T\mathbf{A}+\eta\Pb\Big)^{-1}\mathbf{A}^T\bigg)}     \right\}
\end{equation}
for $\Pb=(\mathbf{I}+(\ell/h)^2\mathbf{L}_{\mbox{\scriptsize2D}})^{\nu+d/2}$.

In practice, $\nu$ is often fixed \cite{khaledi2009empirical,roininen2016hyperpriors} and $\ell$ is either estimated manually or by using the fully Bayesian approach, which involves Markov chain Monte Carlo (MCMC) \cite{robert2013monte} sampling. This requires setting up hyperprior distributions and can be time consuming, subjective and unintuitive, so we present a new method for selecting these hyperparameters next.

\section{The Semivariogram Method for Estimating $\boldsymbol{\nu}$ and $\boldsymbol{\ell}$}
\label{sec:variogram}

In the inverse problem formulation above, the components of the vector $\boldsymbol{x}$ correspond to values of an unknown function $x$ at numerical mesh points within a spatial region $\Omega$. This motivates using methods from spatial statistics to estimate the Whittle-Mat\'ern prior hyperparameters $\nu$ and $\ell$. One such method uses a variogram, and a corresponding semivariogram \cite{schabenberger2017statistical}, which requires the assumption of {\em intrinsic stationarity}, i.e., that the elements of $\xb$ have constant mean and the variance of the difference between the elements is constant throughout the region. This is a weaker assumption than is required by many other parameter estimation tools, which is one of the reasons variograms have become popular in spatial statistical applications \cite{cressie2015statistics}, and it is the reason we use semivariograms here. Although the use of semivariograms for estimating parameters to determine a covariance structure is commonly used in spatial statistics, this is, to our knowledge, the first time these tools have been used to estimate prior hyperparameters for use in inverse problems.

The semivariogram is defined by $\gamma(\boldsymbol{r})=\frac{1}{2}\mbox{Var}[Z(\boldsymbol{u}_i)-Z(\boldsymbol{u}_j)],$ where $\boldsymbol{r}=\boldsymbol{u}_i-\boldsymbol{u}_j$ and $\{Z(\boldsymbol{u}):\boldsymbol{u}\in \Omega\subset\mathbb{R}^d\}$ is a spatial process.
Due to our stationarity assumption, $\mbox{Var}[Z(\boldsymbol{u}_i)]=\mbox{Var}[Z(\boldsymbol{u}_j)]=\sigma^2$, which we use to derive the following alternative expression for $\gamma(\rb)$:
\begin{eqnarray*}
\gamma(\boldsymbol{r})
&=\frac{1}{2}\Big(\mbox{Var}[Z(\boldsymbol{u}_i)]+\mbox{Var}[Z(\boldsymbol{u}_j)]-2\mbox{Cov}[Z(\boldsymbol{u}_i),Z(\boldsymbol{u}_j)]\Big)\\
& = \sigma^2-\mbox{Cov}[Z(\boldsymbol{u}_i),Z(\boldsymbol{u}_j)].
\end{eqnarray*}
Thus, the semivariogram simplifies to the difference between the variance in the region and the covariance between two points with a difference $\boldsymbol{r}$. The variogram is formally defined as $2\gamma(\boldsymbol{r})$, hence the terms variogram and semivariogram are often used interchangeably. To remain consistent, we will continue to refer to $\gamma(\rb)$ as a semivariogram throughout the paper.

We now need a way to estimate the semivariogram from given data. For this, we use what is known as the sample, or empirical, semivariogram. Assuming that $Z(\ub)$ is isotropic, so that $r=\Vert \boldsymbol{r}\Vert=\Vert \boldsymbol{u}_i-\boldsymbol{u}_j\Vert$, then the empirical semivariogram can be expressed
\begin{equation}
\label{empiricalSV}
\hat{\gamma}(r)=\frac{1}{2n(r)}\sum_{(i,j)\left|\Vert\boldsymbol{u}_{i}-\ub_j\Vert=r\right.}\hspace{-6mm}[z(\boldsymbol{u}_i)-z(\boldsymbol{u}_j)]^2,
\end{equation}
where $z(\boldsymbol{u})$ is a realization of $Z(\boldsymbol{u})$, and $n(r)$ is the number of points that are separated by a distance $r$. The $\hat{\gamma}(r)$ values are often referred to as the semivariance values. In a typical semivariogram, the semivariance values increase as $r$ increases since points tend to be less similar the further apart they are, which increases the variance of their differences.

Although the empirical semivariogram is useful in obtaining semivariance values from data, it is not ideal for modeling data for various reasons (see \cite{cressie2015statistics} for details), thus it is typical to fit a semivariogram model to the empirical semivariogram. Since our prior distribution for $\boldsymbol{x}$ has a Mat\'ern covariance, we will use the theoretical Mat\'ern semivariogram model \cite{matern2013spatial, stein2012interpolation} given by
\begin{eqnarray}
\label{maternSV}
\fl\gamma(r,\boldsymbol{\theta})=\cases{0&if $r=0$\\
a_0+(\sigma^2-a_0)\left[1-\frac{1}{2^{\nu-1}\Gamma(\nu)}(r/\ell)^\nu K_\nu( r/\ell)\right] &if $r>0$\\}
\end{eqnarray}
where $a_0\ge0$ is the nugget, $\sigma^2\ge a_0$ is the sill, and $\boldsymbol{\theta}=(a_0,\sigma^2,\nu,\ell)$. The nugget is the term given to the semivariance value at a distance just greater than zero and the sill is the total variance contribution or the semivariance value where the model levels out. The sill, $\sigma^2$, is also the variance parameter in the Mat\'ern covariance function (\ref{equ:matern}). We can estimate $a_0$, $\sigma^2$, $\nu$, and $\ell$ by fitting semivariogram models to the empirical semivariogram.

There are a number of ways to fit the semivariogram model to the empirical semivariogram. We use weighted least squares, as is commonly done \cite{cressie2015statistics}, choosing the $\thetab$ that minimizes
\begin{eqnarray}
\label{wss}
W(\thetab)=\sum_{r}\frac{n(r)}{2[\gamma(r,\boldsymbol{\theta})]^2}[\hat{\gamma}(r)-\gamma(r,\boldsymbol{\theta})]^2.
\end{eqnarray}
To minimize $W(\thetab)$, we adapt the MATLAB codes from \cite{Matlab_Variogram,Matlab_Variogramfit}. More specifically, we adapt \cite{Matlab_Variogram} for computing the empirical semivariance $\hat{\gamma}(r)$ and we adapt \cite{Matlab_Variogramfit} for minimizing $W(\thetab)$. Although it is possible to optimize both $\nu$ and $\ell$ continuously, we will require $\nu+d/2$ to be an integer. Weighted least squares, in general, performs well when finding optimal estimates for $a_0, \sigma^2$, and $\ell$ for given empirical semivariogram values when $\nu$ is fixed, but not when $\nu$ is also free to vary (most software requires a fixed $\nu$ value). To combat this issue, and to ensure $\nu+d/2$ is an integer, we cycle through various fixed values of $\nu$ to obtain estimates for the other parameters and their weighted least squares value. We then choose the $\boldsymbol{\theta}$ with the smallest $W(\boldsymbol{\theta})$.

\begin{figure}
\centering
\raisebox{5.05mm}{\includegraphics[width=2.03in]{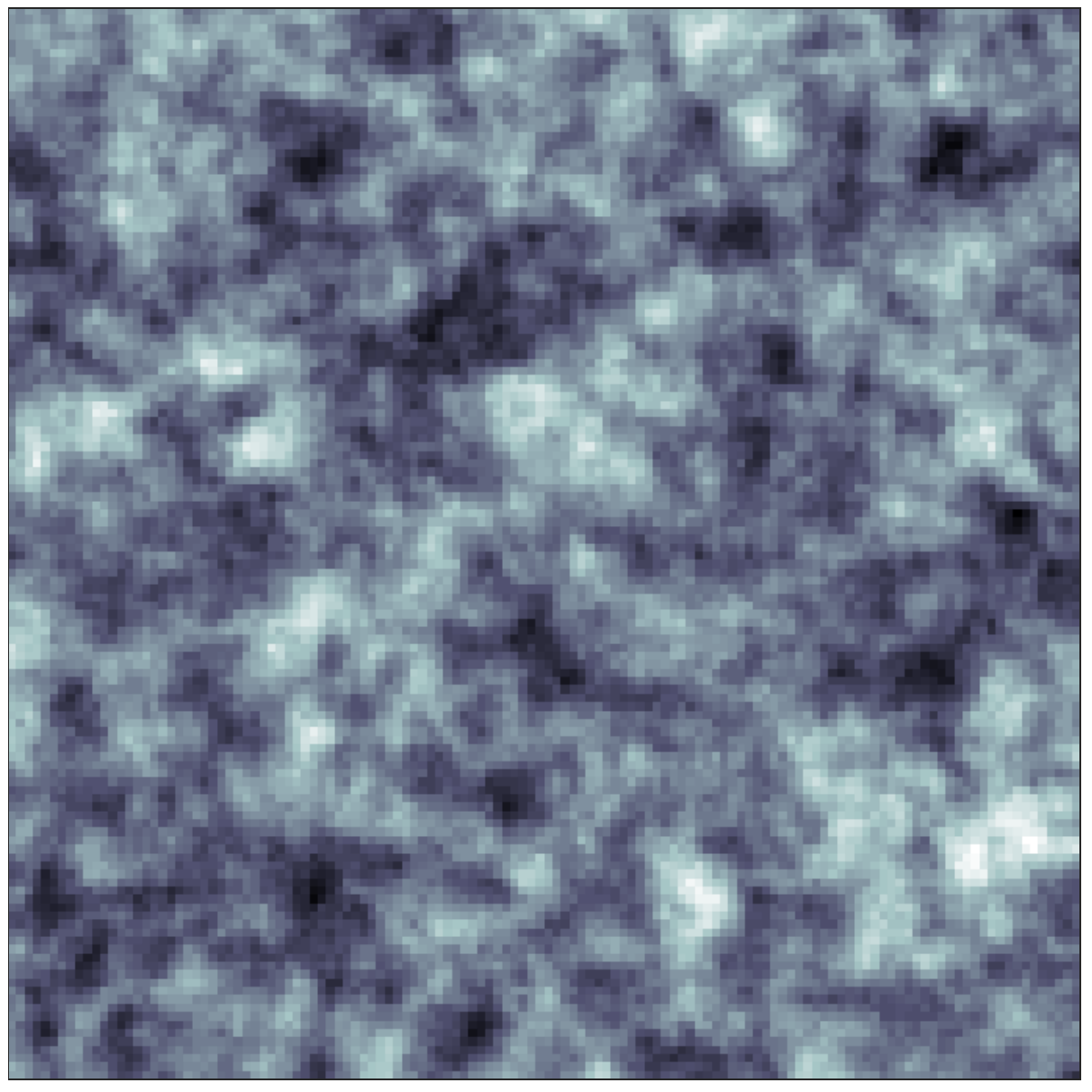}}\quad
\includegraphics[width=2.28in]{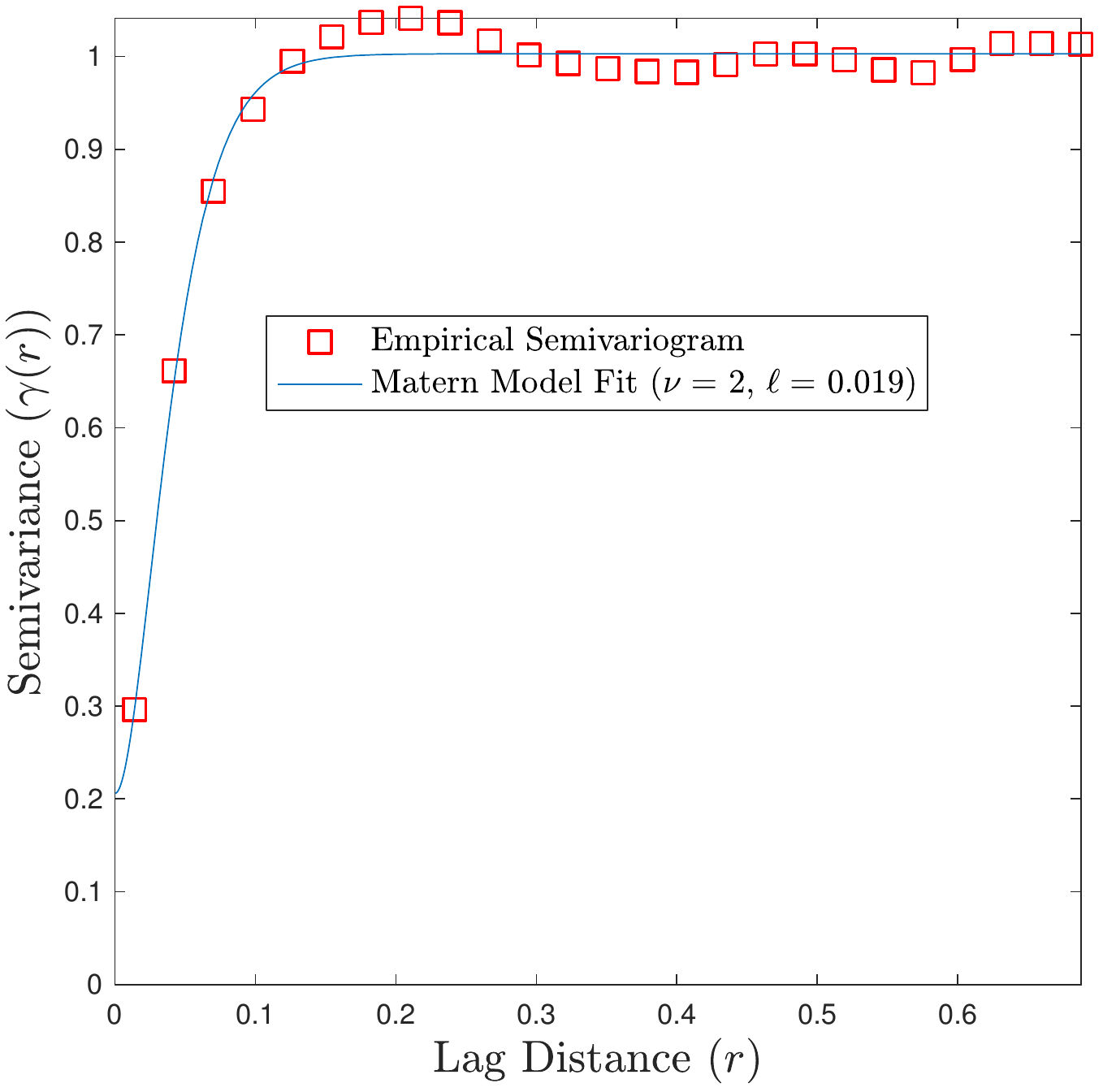}
\caption{Semivariogram. A randomly generated spatial field is shown on the left and the empirical semivariogram, along with the Mat\'ern model fit, is given on the right. The fitted hyperparameters are $\nu=2$ and $\ell=0.019$, which corresponds to a practical range of $0.102$. }
\label{variogram}
\end{figure}

For an illustration, we generated a random field, shown on the left side of Figure \ref{variogram}, and fit a semivariogram to the field. The optimized parameters of the model are $\nu=2$ and $\ell=0.019$, which corresponds to a practical range of $0.102$. Thus, the values of the field are nearly independent a tenth of the way across the region. The sill and nugget are estimated to be $\sigma^2=1.003$ and $a_0=0.206$, respectively. A plot of the resulting fitted Mat\'ern semivariogram model is given on the right side of Figure \ref{variogram}.

The values of $\nu$ and $\ell$ from $\thetab=(a_0,\sigma^2,\nu,\ell)$ obtained by fitting the Mat\'ern semivariogram model to a spatial field, as described in the previous paragraph, can be used to define the Whittle-Mat\'ern prior (\ref{equ:prior}). The sill, $\sigma^2$, and the nugget, $a_0$, are not especially useful outside of fitting the semivariogram model because they do not correspond to any hyperparameter in (\ref{equ:prior}). They are helpful only in determining the best estimates for $\nu$ and $\ell$. Any contribution these parameters may have made to the prior distribution will be accounted for in the regularization parameter, $\alpha$. Therefore, after fitting the semivariogram models, $\sigma^2$ and $a_0$ are discarded.

\begin{algorithm}
\caption{The semivariogram method for MAP estimation with Whittle-Mat\'ern prior.}
\label{alg1}
0. Estimate $\thetab=(a_0,\sigma^2,\nu,\ell)$ by fitting a Mat\'ern semivariogram model to $\boldsymbol{b}$.\\
1. Define the prior (\ref{equ:prior}) using $\nu$ and $\ell$, compute $\alpha$ using (\ref{equ:GCV}), and compute $\boldsymbol{x}_\alpha$ using (\ref{equ:MAP}).\\
2. Update $\thetab=(a_0,\sigma^2,\nu,\ell)$ by fitting a Mat\'ern semivariogram model to $\boldsymbol{x}_\alpha$.\\
3. Return to step 1 and repeat until $\nu$ and $\ell$ stabilize.
\end{algorithm}

With estimates for $\nu$ and $\ell$ in hand, the MAP estimator, $\xb_\alpha$, can then be computed as in Section \ref{sec:MAP}, from which we can recompute $\thetab$ by fitting the Mat\'ern semivariogram model to the empirical semivariogram values of $\xb_\alpha$. Repeating this process iteratively yields Algorithm \ref{alg1}. 
%\noindent{\bf Algorithm:} {\em The Semivariogram Method for MAP Estimation with Whittle-Mat\'ern Prior:}
%\begin{enumerate}
%\item[0.] Estimate $\thetab=(a_0,\sigma^2,\nu,\ell)$ by fitting a Mat\'ern semivariogram model to $\boldsymbol{b}$.
%
%\item[1.] Define the prior (\ref{equ:prior}) using $\nu$ and $\ell$, compute $\alpha$ using (\ref{equ:GCV}), and compute $\boldsymbol{x}_\alpha$ using (\ref{equ:MAP}).
%
%\item[2.] Update $\thetab=(a_0,\sigma^2,\nu,\ell)$ by fitting a Mat\'ern semivariogram model to $\boldsymbol{x}_\alpha$.
%
%\item[3.] Return to step 1 and repeat until $\nu$ and $\ell$ stabilize.
%\end{enumerate}
Recall that $\boldsymbol{b}$ is a vector of measurements, which will usually be noisy or have some missing values, and each element of $\boldsymbol{b}$ has a corresponding spatial position. Since $\nu$ is being optimized discretely to ensure that $\beta=\nu+d/2$ is an integer, convergence will be met when $\nu_j-\nu_{j-1}=0$ where $\nu_j$ is the $\nu$ value fit in the $j$th iteration. Then $\ell$ is said to have converged when $|\ell_j-\ell_{j-1}|/\ell_{j-1}<\varepsilon$ with $\varepsilon$ determined by the user. In this paper, we will consider $\ell$ to have converged when the relative difference is less than $0.01$, which usually takes fewer than three iterations to achieve.

The semivariogram method is essentially a parametric empirical Bayes method \cite{casella1985introduction} for point estimation. We have a distributional assumption on $\xb$, but no prior distributions are assumed for $\nu$ or $\ell$. 
The hyperparameters are instead estimated by iteratively fitting semivariograms to the data. 

\subsection{Numerical Experiments}
We now implement the semivariogram method on a two-dimensional deblurring and inpainting example. Recall that the connection between the Mat\'ern covariance and the Whittle-Mat\'ern prior depends on a stationarity assumption, which the following example may not exhibit. For simplicity, we will still assume stationarity and acknowledge that future work should be done in the case when no stationarity is present. Additionally, the numerical examples given in this paper all use color images. In our analysis, we will assume independence in the color bands and obtain priors and reconstructions for each one individually.

\subsubsection{Results}
\label{iso_results}

\begin{figure}[b!]
\centering
\includegraphics[width=2.5in]{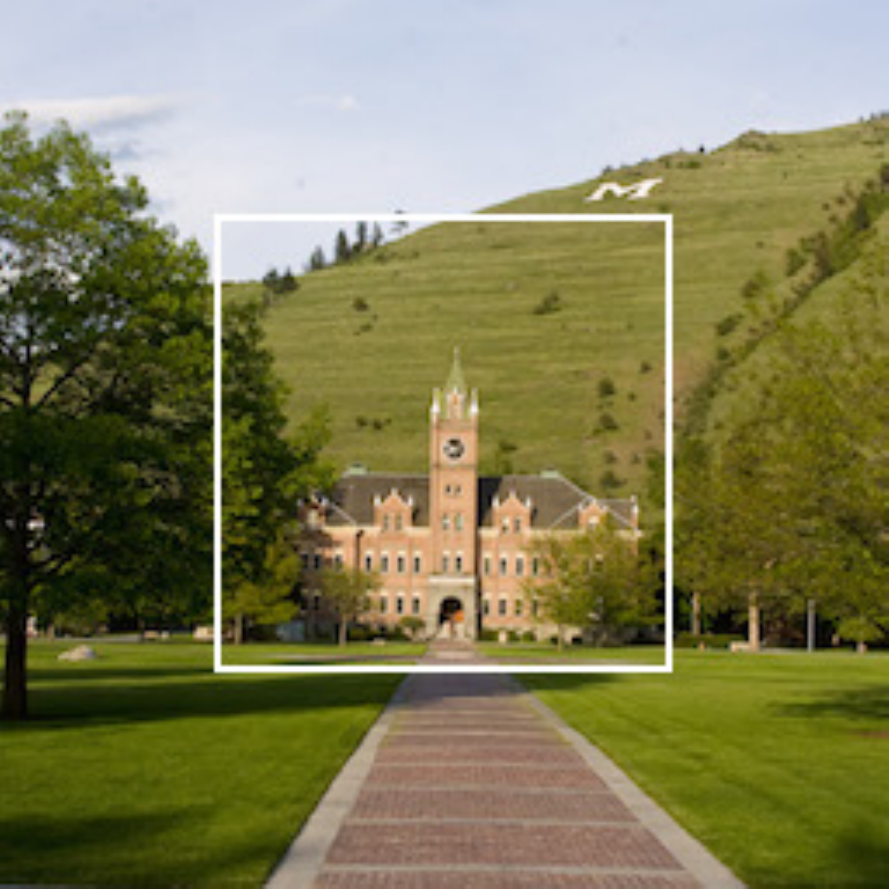}
\caption{Full $256\times 256$ image of Main Hall at the University of Montana with $128\times 128$ subimage.}
\label{2dmainhall}
\end{figure}

In this example, we assume periodic boundary conditions on the extended domain, but due to the restriction from the extended domain $\overline\Omega$ to $\Omega$, circulant structure is lost in the forward model matrix, and hence, linear system solves must be done using an iterative method. As in \cite[Section 3.1.3]{bardsley2018}, we use preconditioned conjugate gradient (PCG) iteration, both for computing $\alpha$ and for computing $\boldsymbol{x}_\alpha$. We attempt to deblur and demask a $128\times128$ image of Main Hall on the University of Montana (UM) campus. To do this, we begin with a $256\times 256$ image, given in Figure \ref{2dmainhall}, and then restrict to the center $128\times128$ image. This smaller image in the middle will be thought of as being on a domain $\Omega=[0,1]\times[0,1]$ and the larger, full image will then be defined on $\overline\Omega=[-0.5,1.5]\times[-0.5,1.5]$.

\begin{figure}
\centering
\includegraphics[width=1.65in]{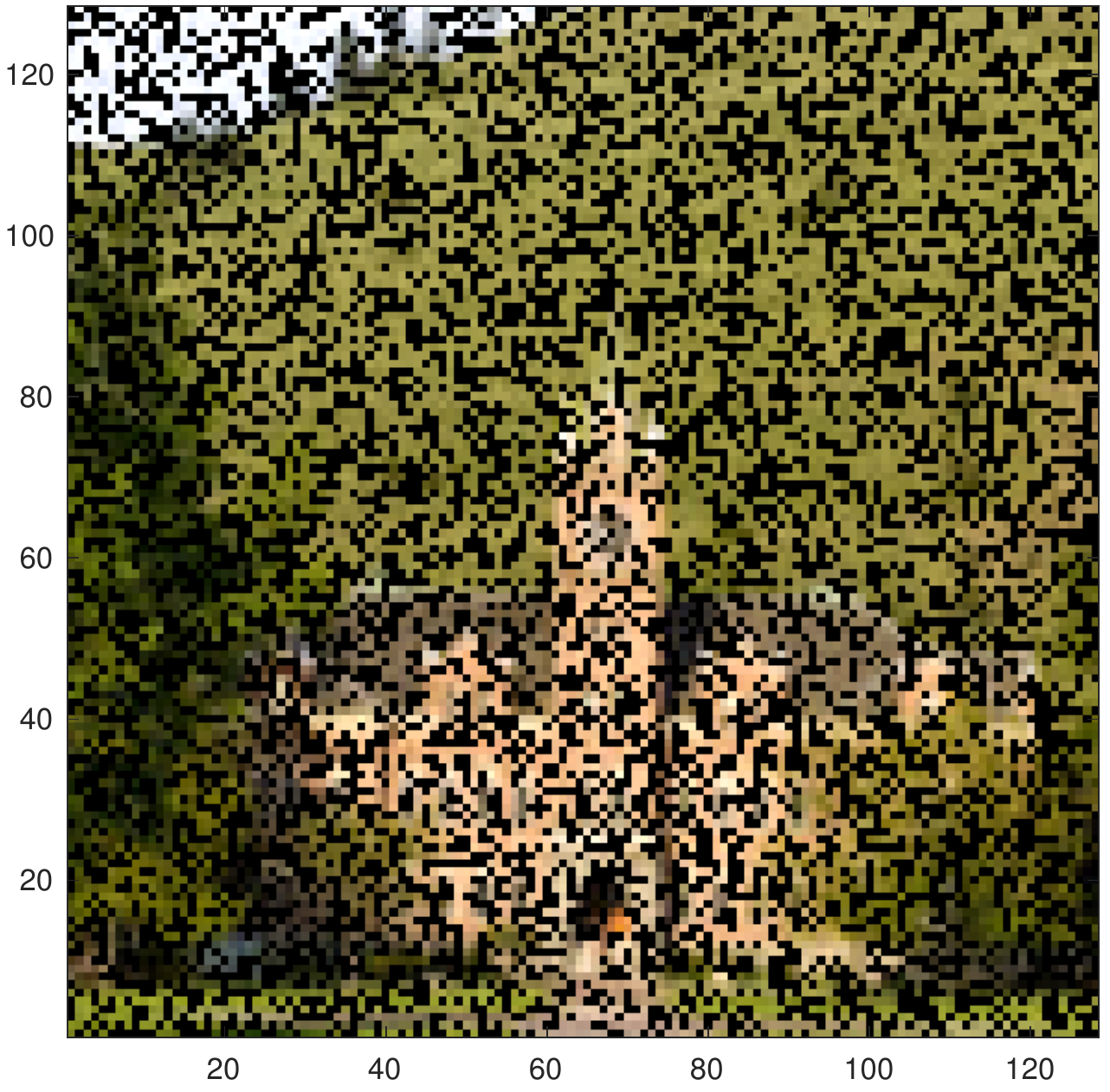}
\includegraphics[width=1.65in]{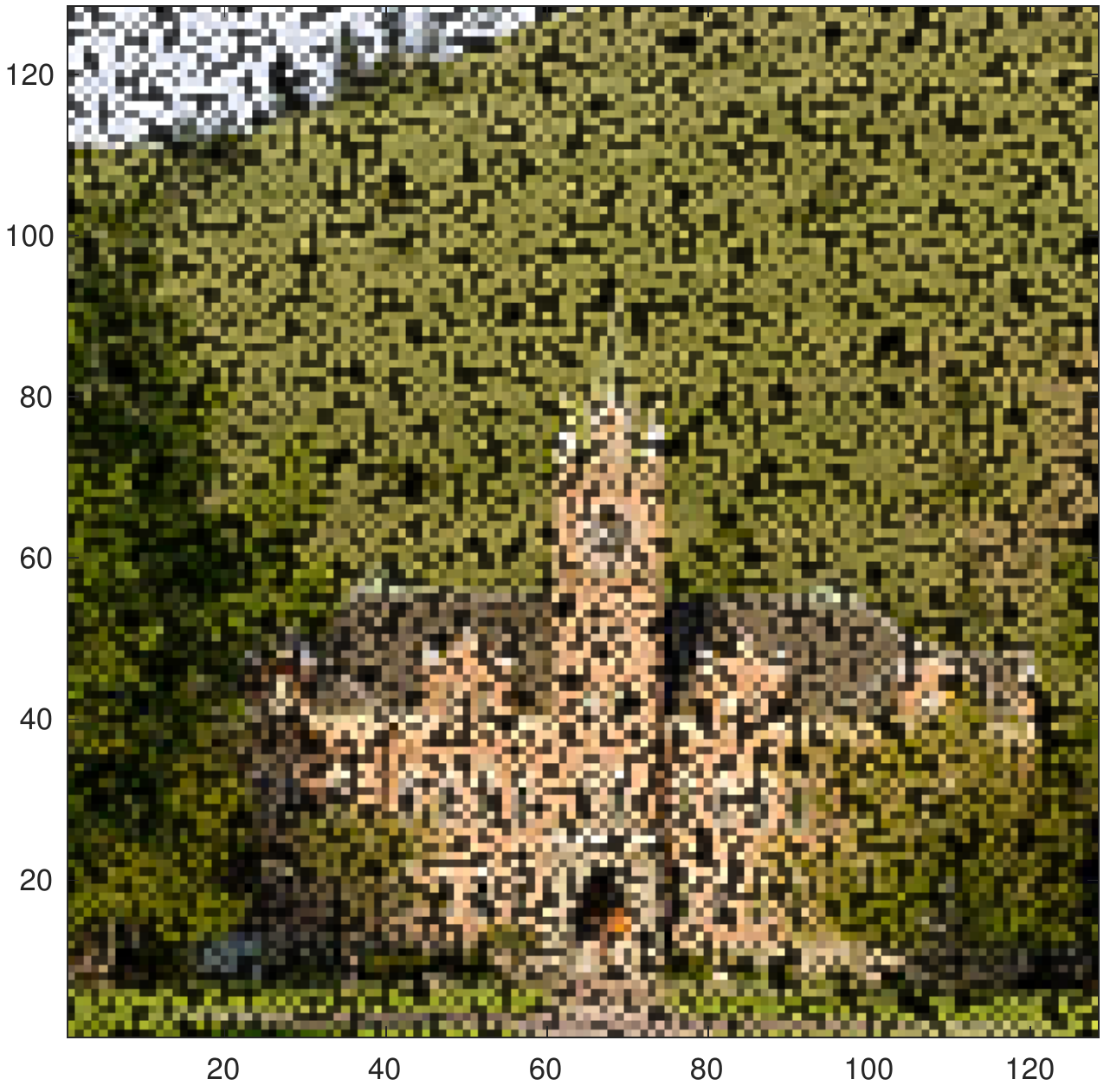}
\includegraphics[width=1.65in]{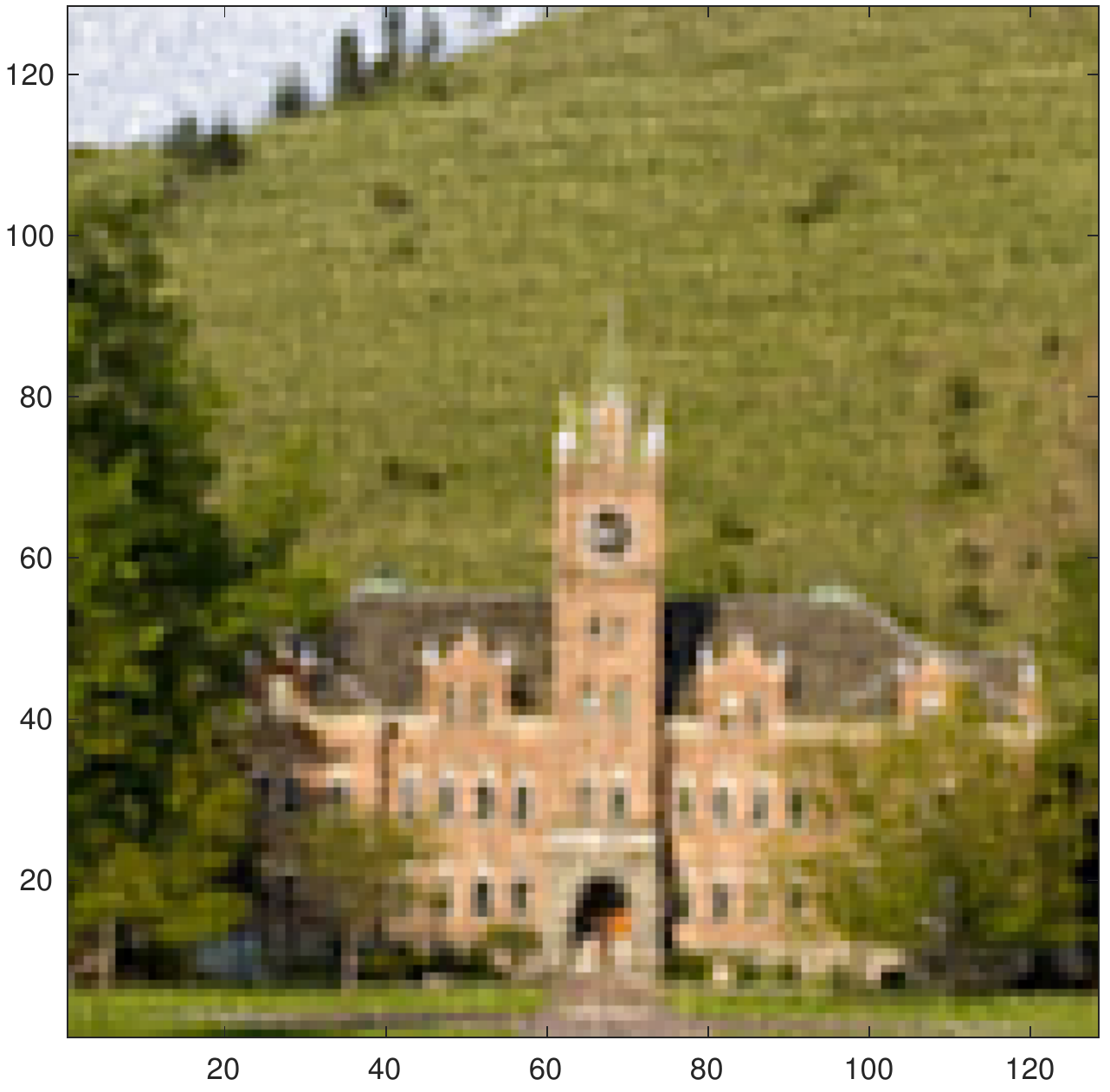}
\caption{Two-dimensional image deblurring test case. On the left is a plot of the blurred, masked, and noisy data; in the middle is a plot of the Tikhonov solution; and on the right is a plot of the solution obtained using the Whittle-Mat\'ern prior with $\nu=1$ and $\ell=0.0364, 0.0313$ and $0.0543$ for red, green and blue intensities, respectively.}
\label{2dexample}
\end{figure}

To obtain $\bb$, we first perform a slight blurring operation on the full 256$\times$256 true image plotted in Figure \ref{2dmainhall}. Since this is a color image, the deblurring process is done individually for the red, green, and blue intensity arrays. We then restrict to the central $128\times128$ pixels (with boundaries denoted in Figure \ref{2dmainhall}) and randomly remove 40\% of the pixels to obtain the masked, and moderately blurry image on the left in Figure \ref{2dexample}.

We seek an estimate of $\xb$ in the same central subregion. Omnidirectional semivariograms with 25 approximately equally spaced grid points in $0 <r<\sqrt{2}/10$ are used. We chose $\sqrt{2}/10$ as a cutoff because it balances the need to capture the covariance structure at short distances, which are well-known to be the most important \cite{cressie2015statistics}, with those at longer distances. When fitting semivariograms to the masked image, the removed entries will not be considered or else the correlation would be strongly influenced by those entries.

The semivariogram method is used to obtain $\nu=1$ for each color band, $\ell=0.0364, 0.0313$ and $0.0543$ for the red, green, and blue intensities, respectively, and $\alpha=0.0023, 0.0018$ and $5.47\times 10^{-5}$. Convergence was met in two iterations for each color intensity. We also computed the Tikhonov solution, as defined in \cite[Section 3.1.3]{bardsley2018}, for which the prior covariance is equal to a scalar multiple of the identity matrix. The Tikhonov $\alpha$ values for all three color bands were around $0.0004$. Note that for both of these reconstructions, the regularization parameter, $\alpha$, was optimized using the highest correlation between the solution and the true image rather than chosen by GCV to ensure that any differences in the solutions is due to the method and not a poorly-chosen regularization parameter.

The two solutions are plotted in Figure \ref{2dexample}. It is clear that the solution that used the Whittle-Mat\'ern prior is the superior reconstruction. The correlation between $\xb_\alpha$ and $\xb$, the true image, is 0.950. While the Tikhonov solution is able to remove the blur, it performs inpainting poorly since each pixel value is assumed independent of one another due to the identity covariance matrix.

\subsection{Discussion}
Compared to the fully Bayesian method, the semivariogram procedure has some key advantages. This technique produces competitive solutions and clearer interpretations of the hyperparameters $\nu$ and $\ell$, and it can inform how far to extend the domain to maintain a connection with the Mat\'ern covariance. Additionally, the computation time is only a fraction of what is needed to compute an adequate number of MCMC samples. In our implementation of the example above, the semivariogram method was more than 20 times faster than the fully Bayesian MCMC method. Finally, it is not trivial to sample from a complex model such as this one without significant autocorrelation, whereas sampling is not needed for the semivariogram method. 

The primary disadvantage is that we lose uncertainty quantification. We also have to calculate $\alpha$, the regularization parameter, using other techniques like GCV. One other shortcoming to the semivariogram method, as described in this section, is the fact that it requires the field or image to be isotropic. In the next section, we extend these results to anisotropic fields.

\section{Geometric Anisotropy}
\label{sec:anis}

The solution to (\ref{equ:SPDE}) is an isotropic Gaussian field, which means the correlation length is the same in every direction. This isotropy assumption is often not satisfied and so it will be useful to have an alternate SPDE formulation for the case when correlation lengths differ with direction. This is known as geometric anisotropy \cite{cressie2015statistics}. The groundwork for constructing priors that can model anisotropy has been laid in works such as \cite{roininen2014whittle,lindgren2011explicit,hale2013implementing}.

\subsection{Anisotropic SPDE}

We will derive an anisotropic SPDE that can be used in similar way that (\ref{equ:SPDE}) was used in the prior modeling in the isotropic case. We will only consider the two-dimensional case, but results can be extended to $d>1$ dimensions. In two dimensions, for a Gaussian field with correlation length $\ell_1$ in the direction of the angle $\theta$, where $-\pi/2<\theta\le\pi/2$ is measured counter-clockwise from the $x$-axis, and correlation length $\ell_2$ in the direction perpendicular to $\theta$, we can make the following change of variables from isotropic to anisotropic coordinates:
\begin{equation}
\label{iso_coord_transform}
\boldsymbol{w}=\left[\begin{array}{cc}\cos\theta & -\ell_2/\ell_1\sin\theta\\
\sin\theta & \phantom{-}\ell_2/\ell_1\cos\theta\end{array}\right]
\left[\begin{array}{@{}c@{}}
u_1\\ u_2
\end{array}\right]
\end{equation}
and thus 
\begin{eqnarray*}
&w_1(u_1,u_2)=\cos\theta u_1-\ell_2/\ell_1\sin\theta u_2\\
&w_2(u_1,u_2)=\sin\theta u_1+\ell_2/\ell_1\cos\theta u_2.
\end{eqnarray*}

We will apply the change of variables (\ref{iso_coord_transform}) to both sides of (\ref{equ:SPDE}) to obtain the analogous anisotropic SPDE. The Laplacian on the left-hand side can be altered using the chain rule:
\begin{eqnarray*}
\frac{\partial}{\partial u_1}&=\frac{\partial}{\partial w_1}\frac{\partial w_1}{\partial u_1}+\frac{\partial}{\partial w_2}\frac{\partial w_2}{\partial u_1} \quad \mbox{and}\\
\frac{\partial}{\partial u_2}&=\frac{\partial}{\partial w_1}\frac{\partial w_1}{\partial u_2}+\frac{\partial}{\partial w_2}\frac{\partial w_2}{\partial u_2},
\end{eqnarray*}
which means 
\begin{eqnarray*}
\fl
\frac{\partial^2}{\partial u_1^2}&=\left(\frac{\partial^2}{\partial w_1^2}\frac{\partial w_1}{\partial u_1}+\frac{\partial^2}{\partial w_1\partial w_2}\frac{\partial w_2}{\partial u_1}\right)\frac{\partial w_1}{\partial u_1}+\left(\frac{\partial^2}{\partial w_1\partial w_2}\frac{\partial w_1}{\partial u_1}+\frac{\partial^2}{\partial w_2^2}\frac{\partial w_2}{\partial u_1}\right)\frac{\partial w_2}{\partial u_1}\\
\fl
&=\cos^2\theta\frac{\partial^2}{\partial w_1^2}+2\sin\theta\cos\theta\frac{\partial^2}{\partial w_1\partial w_2}+\sin^2\theta\frac{\partial^2}{\partial w_2^2}
\end{eqnarray*}
and
\begin{eqnarray*}
\fl
\frac{\partial^2}{\partial u_2^2}&=\left(\frac{\partial^2}{\partial w_1^2}\frac{\partial w_1}{\partial u_2}+\frac{\partial^2}{\partial w_2\partial w_1}\frac{\partial w_2}{\partial u_2}\right)\frac{\partial w_1}{\partial u_2}+\left(\frac{\partial^2}{\partial w_1\partial w_2}\frac{\partial w_1}{\partial u_2}+\frac{\partial^2}{\partial w_2^2}\frac{\partial w_2}{\partial u_2}\right)\frac{\partial w_2}{\partial u_2}\\
\fl
&=(\ell_2/\ell_1)^2\sin^2\theta\frac{\partial^2}{\partial w_1^2}-2(\ell_2/\ell_1)^2\sin\theta\cos\theta\frac{\partial^2}{\partial w_1\partial w_2}+(\ell_2/\ell_1)^2\cos^2\theta\frac{\partial^2}{\partial w_2^2}.
\end{eqnarray*}

The right hand side of (\ref{equ:SPDE}) is updated by changing the coordinates of the white noise. The inverse transformation of (\ref{iso_coord_transform}) is
\begin{equation}
\label{anis_coord_transform}
\boldsymbol{u}=\left[\begin{array}{cc}
\cos\theta & \sin\theta\\
-\tau\sin\theta & \tau\cos\theta
\end{array}\right]\left[\begin{array}{@{}c@{}}
w_1\\
w_2
\end{array}\right]:=f(\wb),
\end{equation}
where $\tau=\ell_1/\ell_2$. Now, we define the transformed white noise basis functions as 
\begin{eqnarray*}
\tilde{\phi}_j(\wb)=\psi_j(f(\wb))|\det(J_f(\wb))|^{1/2}=\psi_j(f(\wb))(\ell_1/\ell_2)^{1/2},
\end{eqnarray*}
where $\det(J_f(\wb))$ denotes the determinant of the Jacobian of the transformation $f(\wb)$, which is $(\ell_1/\ell_2)^{1/2}$ in our case. This will preserve the orthonormal properties of the basis functions. Then, appealing to (\ref{equ:white}), 
\begin{eqnarray*}
\mathcal{W}(\boldsymbol{w})&=\sum_{j=1}^\infty\xi_j\tilde{\phi}_j(\boldsymbol{w}),\quad \xi_j\stackrel{iid}{\sim}\mathcal{N}(0,1)\\
&=\sum_{j=1}^\infty\xi_j\phi_j(\ub)(\ell_1/\ell_2)^{1/2}=(\ell_1/\ell_2)^{1/2}\mathcal{W}(\boldsymbol{u}),
\end{eqnarray*}
which means $\mathcal{W}(\boldsymbol{u}) = (\ell_2/\ell_1)^{1/2}\mathcal{W}(\boldsymbol{w})$.

So, taking $\ell=\ell_1$ and making the appropriate substitutions, (\ref{equ:SPDE}) is converted to the anisotropic SPDE:
\begin{equation*}
\fl
\Bigg(1-\hspace{-1mm}\left[(a_\theta^2+b_\theta^2)\scalebox{1.2}{$\frac{\partial^2}{\partial w_1^2}$}+(c_\theta^2+d_\theta^2)\scalebox{1.2}{$\frac{\partial^2}{\partial w_2^2}$}-2(a_\theta c_\theta-b_\theta d_\theta)\scalebox{1.2}{$\frac{\partial^2}{\partial w_1\partial w_2}$}\right]\hspace{-1mm}\Bigg)^{\beta/2}\hspace{-4mm}x(\boldsymbol{w})= (\ell_2/\ell_1)^{1/2}\mathcal{W}(\boldsymbol{w})
\end{equation*}
where $a_\theta=\ell_2\sin\theta$, $b_\theta=\ell_1\cos\theta$, $c_\theta=\ell_2\cos\theta$, and $d_\theta=\ell_1\sin\theta$. For 
\begin{eqnarray*}
\Rb=\left[\begin{array}{cc}
\phantom{-}\ell_1\cos\theta & \ell_1\sin\theta\\
-\ell_2\sin\theta & \ell_2\cos\theta
\end{array}\right],
\end{eqnarray*}
the above SPDE can be written 
\begin{eqnarray}
\label{SPDE_anis}
\left(1-\nabla\cdot \Rb^T\Rb\nabla\right)^{\beta/2}x(\boldsymbol{w})= (\ell_2/\ell_1)^{1/2}\mathcal{W}(\boldsymbol{w}).
\end{eqnarray}
Notice that if $\ell_1=\ell_2$, this SPDE is equivalent to (\ref{equ:SPDE}) with $\ell=\ell_1$.

\subsection{The Gaussian Field Solution of the SPDE (\ref{SPDE_anis})}

Like in the isotropic case, we are interested in the the properties of the solution of (\ref{SPDE_anis}), especially its covariance function. First, we define the anisotropic Mat\'ern covariance function \cite{haskard2007anisotropic} as 
\begin{equation}
\label{anis_matern}
\fl C(r_w)=\sigma^2 \frac{(r_w/\zeta)^{\nu}K_{\nu}(r_w/\zeta)}{2^{\nu-1}\Gamma(\nu)},\ \mbox{with}\ \zeta=\frac{\ell_1}{\sqrt{\cos^2(\psi-\theta)+(\ell_1/\ell_2)^2\sin^2(\psi-\theta)}},
\end{equation}
where $r_w=\Vert \wb_i-\wb_j\Vert$ is the distance between the anisotropic coordinates, $\zeta$ is the new range parameter in the direction of $\psi$, $\ell_1$ is the correlation length in the direction of $\theta$ and $\ell_2$ is the correlation length in the direction perpendicular to $\theta$. Notice that the smoothness parameter, $\nu$, is unaffected.

The remainder of this subsection contains results used to prove the following theorem.
\begin{theorem}
\label{thm:connection_anis}
The solution $x(\boldsymbol{w})$ of (\ref{SPDE_anis}) is a Gaussian field with mean zero and anisotropic Mat\'ern covariance function defined by (\ref{anis_matern}).
\end{theorem}

\noindent {\em Proof.} First, we derive the Green's function for (\ref{SPDE_anis}), which is the solution of
\begin{equation}
\label{eq:anisPDE=dirac}
\left(1-\nabla\cdot \Rb^T\Rb\nabla\right)^{\beta/2}g(\boldsymbol{w},\boldsymbol{v}) = \delta_f(\boldsymbol{v}-\boldsymbol{w}).
\end{equation}
Using (\ref{equ:green}), the solution to (\ref{SPDE_anis}) is given by
\begin{equation}
\label{equ:anisSPDEsolution}
x(\boldsymbol{w})= (\ell_2/\ell_1)^{1/2}\int_{\mathbb{R}^2} g(\boldsymbol{w},\boldsymbol{v})\mathcal{W}(\boldsymbol{v})d\boldsymbol{v},
\end{equation}
which makes $x(\boldsymbol{w})$ a Gaussian field since it is a linear transformation of Gaussian white noise. Be aware that we are still assuming stationarity in our field. To derive the Green's function $g$ in (\ref{equ:anisSPDEsolution}), we first define $g(\boldsymbol{w}):=g(\boldsymbol{w},\boldsymbol{0})$. Then (\ref{eq:anisPDE=dirac}) implies
\begin{equation}
\label{equ:anisSPDEgreen}
\left(1-\nabla\cdot \Rb^T\Rb\nabla\right)^{\beta/2}g(\boldsymbol{w}) = \delta_f(\boldsymbol{w}).
\end{equation}

We would like to change from the anisotropic coordinates $\wb$ to anisotropic coordinates $\ub$ in (\ref{equ:anisSPDEgreen}) so we can use the results from Section \ref{sec:iso_cov}. We again use (\ref{anis_coord_transform}) for the coordinate change and, in a similar fashion as was done earlier, we apply the chain rule to replace 
$\partial^2/\partial w_1^2$, $\partial^2/\partial w_2^2$, and $\partial^2/(\partial w_1\partial w_2)$ 
in $\left(1-\nabla\cdot \Rb^T\Rb\nabla\right)^{\beta/2}$ with partial derivatives in terms of $\ub$. When making this change, the coefficients of $\partial^2/\partial u_1^2$, $\partial^2/\partial u_2^2$, and $\partial^2/(\partial u_1\partial u_2)$ 
are $\ell_1^2$, $\ell_1^2$, and $0$, respectively and so we have $(1-\nabla\cdot \Rb^T\Rb\nabla)g(\boldsymbol{w})=(1-\ell_1^2\Delta)g(\boldsymbol{u})$. Additionally, we can change variables in the Delta function on the right side of (\ref{equ:anisSPDEgreen}) by multiplying by the determinant of the Jacobian of (\ref{anis_coord_transform}): $\ell_1/\ell_2$. Thus, the change of variables transforms (\ref{equ:anisSPDEgreen}) into the equation
\begin{equation}
\label{equ:transformedSPDEgreen}
\left(1-\ell_1^2\Delta\right)^{\beta/2}g(\boldsymbol{u}) = (\ell_1/\ell_2) \delta_f(\boldsymbol{u}),
\end{equation}
which is equivalent to (\ref{equ:SPDEgreen}) up to a constant. Hence, we can apply the results of Section \ref{sec:iso_cov}. Namely, after changing variables, the solution of (\ref{SPDE_anis}) is a Gaussian field with mean zero and the isotropic Mat\'ern covariance function defined by (\ref{equ:matern}). Notice that the constant that multiplies the Delta function on the right-hand side of (\ref{equ:transformedSPDEgreen}) and the constant that multiplies the integral in (\ref{equ:anisSPDEsolution}) will cancel when going through the process of deriving the covariance function since the constant in (\ref{equ:anisSPDEsolution}) gets squared.

We must now make one final change of variables back to $\wb$ from $\ub$ so our covariance function will be in terms of the anisotropic coordinates rather than the isotropic ones. Since the input to the Mat\'ern correlation function must be a distance between isotropic spatial locations, we need to represent an isotropic distance, $r_u$, in terms of the anisotropic coordinates. Consider $\rb:=\wb_i-\wb_j$. Then, defining $r_w:=\Vert\rb\Vert=\Vert\wb_i-\wb_j\Vert$,
\begin{eqnarray*}
\fl r_u&:=\Vert \ub_i-\ub_j\Vert=\left\Vert \left[\begin{array}{cc}
\cos\theta & \sin\theta\\
-\tau\sin\theta & \tau\cos\theta
\end{array}\right](\wb_i-\wb_j)\right\Vert=\left\Vert \left[\begin{array}{cc}
\cos\theta & \sin\theta\\
-\tau\sin\theta & \tau\cos\theta
\end{array}\right]\rb\right\Vert\\
\fl &=\left\Vert \left[\begin{array}{cc}
\cos\theta & \sin\theta\\
-\tau\sin\theta & \tau\cos\theta
\end{array}\right] \left[\begin{array}{cc}
r_1\\ r_2
\end{array}\right]\right\Vert=\left\Vert \left[\begin{array}{cc}
r_1\cos\theta + r_2\sin\theta\\
-r_1\tau\sin\theta + r_2\tau\cos\theta
\end{array}\right]\right\Vert.
\end{eqnarray*}
Now we convert to polar coordinates with $r_1=r_w\cos\psi$ and $r_2=r_w\sin\psi$. Then
\begin{eqnarray*}
r_u&=\left\Vert \left[\begin{array}{cc}
r_w\cos\psi\cos\theta + r_w\sin\psi\sin\theta\\
-r_w\tau\cos\psi\sin\theta + r_w\tau\sin\psi\cos\theta
\end{array}\right]\right\Vert
=\left\Vert \left[\begin{array}{cc}
r_w\cos(\psi-\theta)\\
r_w\tau\sin(\psi-\theta)\\
\end{array}\right]\right\Vert\\
&=r_w\left[\cos^2(\psi-\theta)+\tau^2\sin^2(\psi-\theta)\right]^{1/2}.
\end{eqnarray*}
Therefore, we need to adjust the distance between the vectors $\wb_i$ and $\wb_j$ by $[\cos^2(\psi-\theta)+\tau^2\sin^2(\psi-\theta)]^{1/2}$ in order to get the distances to plug into the isotropic Mat\'ern correlation function. Thus, the isotropic Mat\'ern covariance function has been generalized to the anisotropic case using the same change of variables as in (\ref{iso_coord_transform}).
Adjusting the anisotropic distances is equivalent to defining the anisotropic Mat\'ern covariance function as 
we have in (\ref{anis_matern}).  \hfill $\square$

\subsection{Anisotropic Prior Modeling}

To obtain a sparse representation of the precision matrix for the anisotropic Mat\'ern covariance, we can discretize (\ref{SPDE_anis}) using the standard finite-difference approximations with appropriate boundary conditions. Taking a step size of $h=1/n$  on a uniform mesh, so that $N=n^2$ in two dimensions, yields
\begin{eqnarray*}
\Big[\mathbf{I}+&\frac{1}{h^2}(a_\theta^2+b_\theta^2)(\mathbf{L}\otimes \mathbf{I})+\frac{1}{h^2}(c_\theta^2+d_\theta^2)(\mathbf{I}\otimes\mathbf{L})\\
&-\frac{2}{4h^2}(a_\theta c_\theta-b_\theta d_\theta)(\mathbf{K}\otimes\mathbf{K})\Big]^{\beta/2}\xb=\delta^{-1/2}\boldsymbol{\xi}, \quad\boldsymbol{\xi}\sim \mathcal{N}(\boldsymbol{0},\boldsymbol{I}_N).
\end{eqnarray*}
where $\otimes$ denotes Kronecker product \cite{bardsley2018}. Note that the constant multiplying the white noise term gets absorbed into the $\delta$ hyperparameter.

In the zero boundary condition case,
\begin{eqnarray*}
\fl
\mathbf{L}=\left[\begin{array}{ccccc}
2 & -1 & 0 & \dots & 0 \\
-1 & 2 & -1 & \dots & 0\\
0 & -1 & 2 & \ddots & \vdots \\
\vdots & \ddots & \ddots & \ddots & -1\\
0 & 0 & \dots & -1 & 2 
\end{array}\right]_{n\times n}
\mbox{ and }\quad
\mathbf{K}=\left[\begin{array}{ccccc} 
0 & 1 & 0 & \dots & 0 \\
-1 & 0 & 1 & \dots & 0\\
0 & -1 & 0 & \ddots & \vdots \\
\vdots & \ddots & \ddots & \ddots & 1\\
0 & 0 & \dots & -1 & 0 
\end{array}\right]_{n\times n},
\end{eqnarray*}
and when using periodic boundary conditions, we let $\mathbf{L}(1,n)=\mathbf{L}(n,1)=\mathbf{K}(1,n)=-1$ and $\mathbf{K}(n,1)=1$. Then
\begin{equation}
\label{disc_anis}
\xb\sim \mathcal{N}\left(\boldsymbol{0},\delta^{-1}\mathbf{P}^{-1} \right)
\end{equation}
where
\begin{eqnarray}
\label{equ:P}
\fl \mathbf{P}=\left[\mathbf{I}+\frac{1}{h^2}(a_\theta^2+b_\theta^2)(\mathbf{L}\otimes \mathbf{I})+\frac{1}{h^2}(c_\theta^2+d_\theta^2)(\mathbf{I}\otimes\mathbf{L})-\frac{2}{4h^2}(a_\theta c_\theta-b_\theta d_\theta)(\mathbf{K}\otimes\mathbf{K})\right]^{\beta}.
\end{eqnarray}
In order to retain sparsity in $\mathbf{P}$, we will again require that $\beta=\nu+d/2$ be an integer. Additionally, like in the isotropic case, an extension of the computational domain is required to maintain a connection between (\ref{equ:P}) and (\ref{anis_matern}).

Now that we have a prior covariance matrix that maintains a connection to the anisotropic Mat\'ern covariance, we return to the MAP estimator, which can be computed by solving
\begin{eqnarray}
\boldsymbol{x}_\alpha
&=\mbox{arg}\min_{\boldsymbol{x}}\left\{\frac{1}{2}\Vert \mathbf{A}\boldsymbol{x}-\boldsymbol{b}\Vert^2 +\frac{\alpha}{2}\boldsymbol{x}^T\mathbf{P}\boldsymbol{x}\right\}\nonumber\\
&=\left(\mathbf{A}^T\mathbf{A}+\alpha\mathbf{P}\right)^{-1}\mathbf{A}^T\boldsymbol{b},\label{equ:MAP2}
\end{eqnarray}
where $\alpha=\delta/\lambda$ and $\mathbf{P}$ is as in (\ref{equ:P}).

\subsection{Directional Semivariograms}

When fitting semivariograms to a spatial field, intrinsic stationarity and isotropy is assumed. In our case, we are still assuming intrinsic stationarity, but our field is anisotropic. Thus, a change must be made to our field before fitting a semivariogram to obtain an estimate for $\ell_1$ and $\ell_2$. We again use (\ref{anis_coord_transform}), the inverse of the change used in (\ref{iso_coord_transform}).
Using the same argument that was used when transforming the Green's function PDE from anisotropic coordinates in (\ref{equ:anisSPDEgreen}) to isotropic coordinates in (\ref{equ:transformedSPDEgreen}), it is not difficult to show that (\ref{SPDE_anis}) is transformed to
\begin{eqnarray*}
(1-\ell_1^{2}\Delta)^{(\nu+d/2)/2}x(\boldsymbol{u})=\mathcal{W}(\ub),
\end{eqnarray*}
which is equivalent to (\ref{equ:SPDE}) with $\ell=\ell_1$.

\begin{figure}
\centering
\includegraphics[width=2.25in]{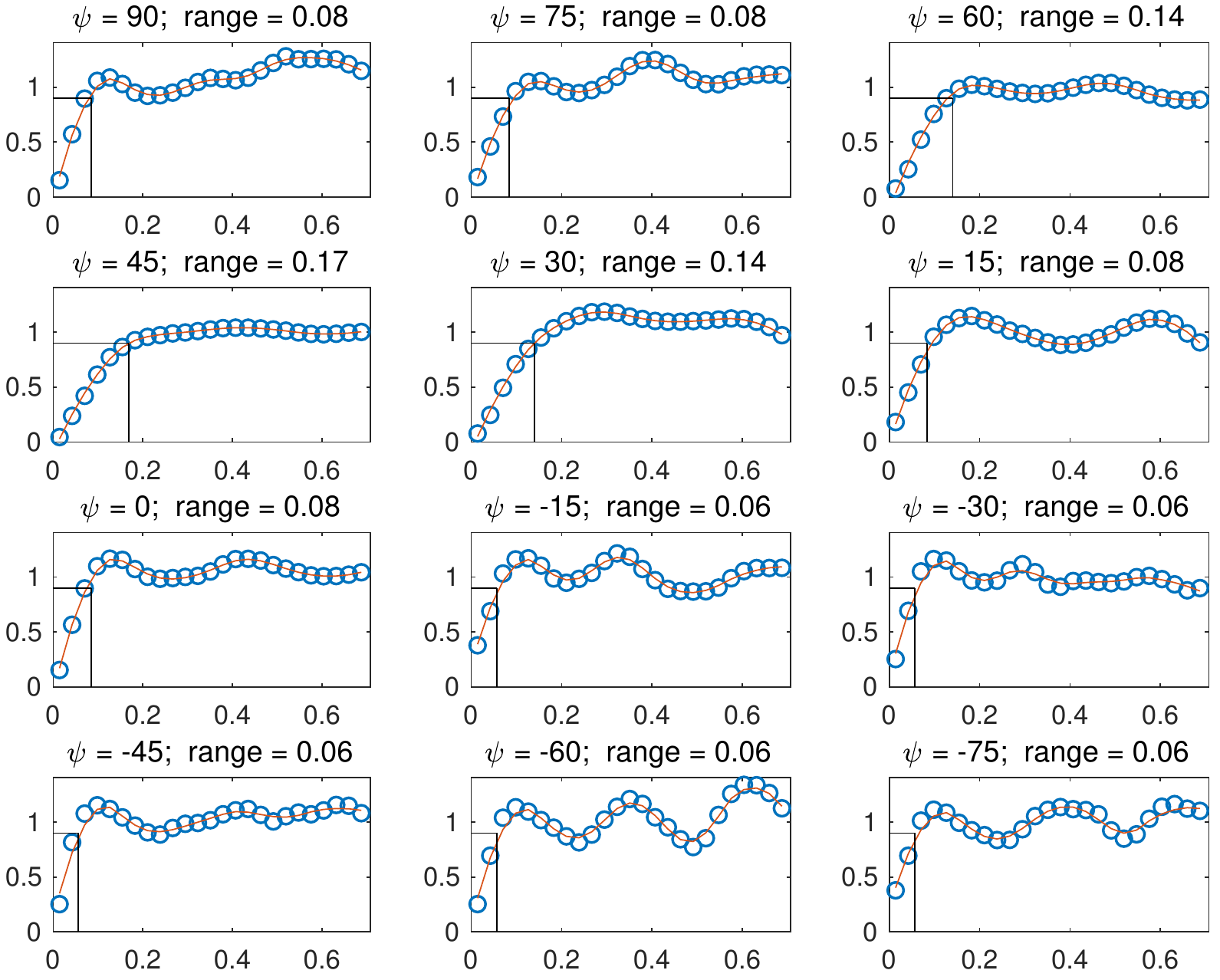}\quad
\raisebox{1.5mm}{\includegraphics[width=1.3in]{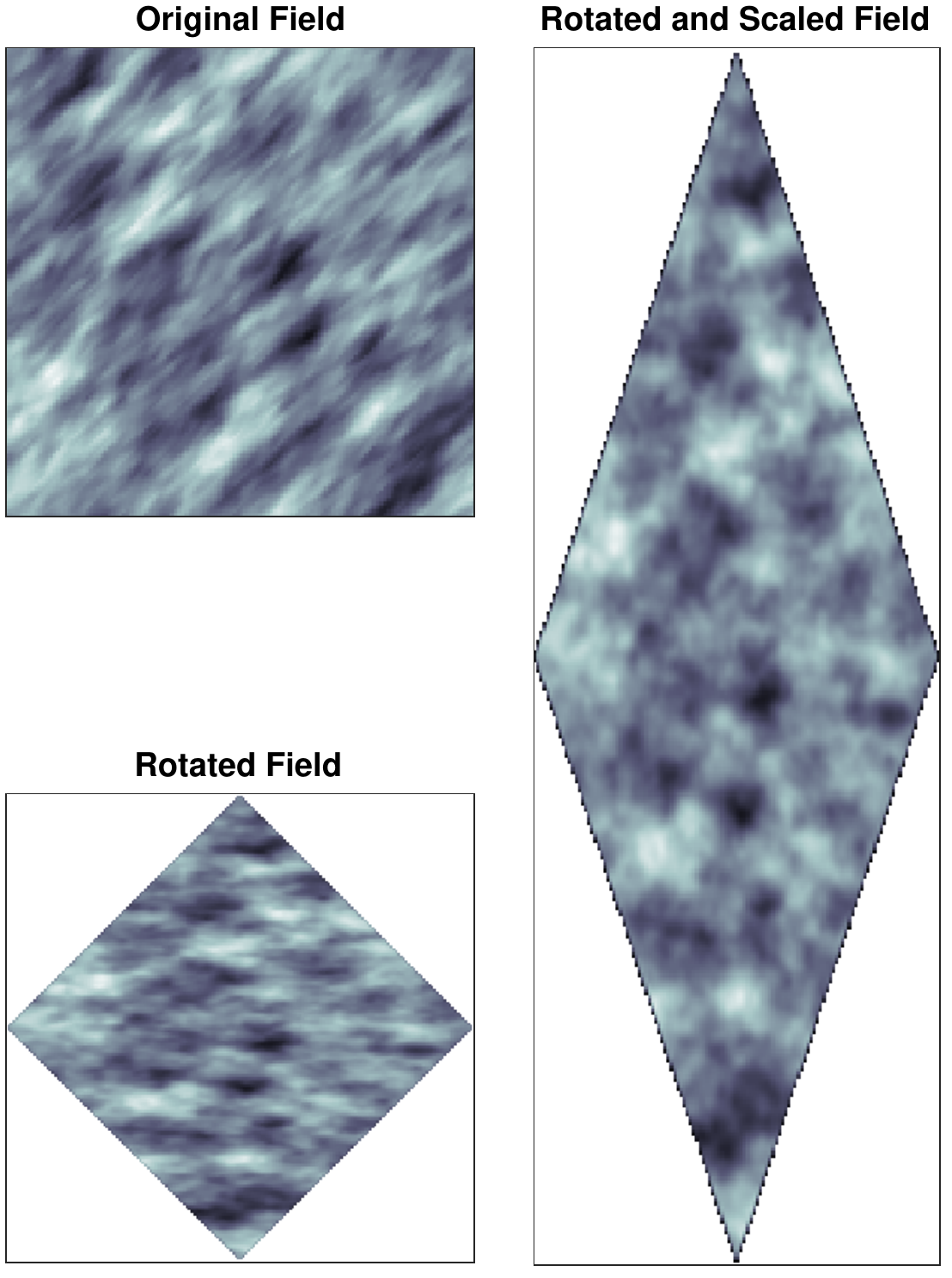}}\quad
\includegraphics[width=2.25in]{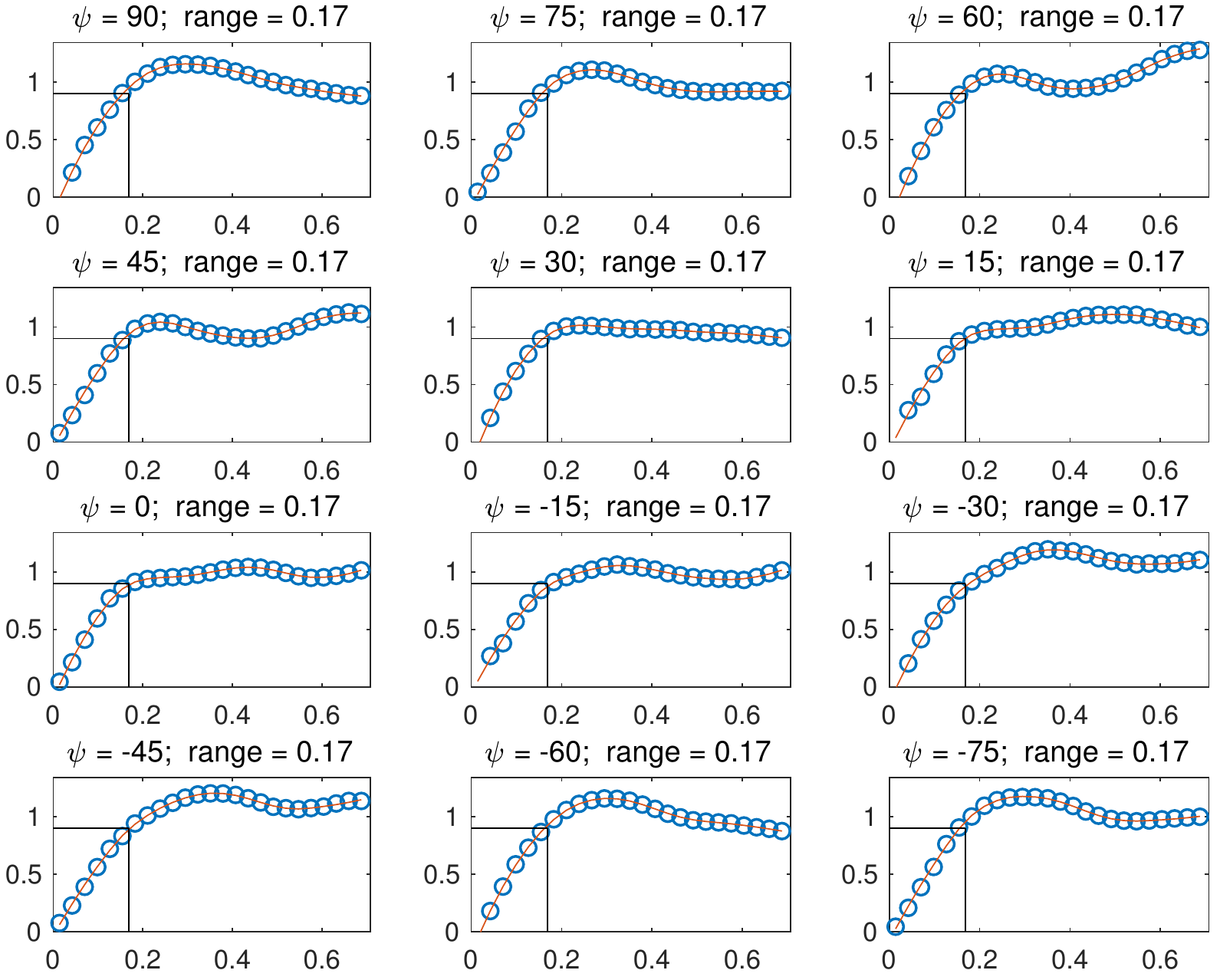}
\caption{Directional semivariograms. The directional semivariograms for the original, anisotropic field are shown on the left. For each of the 12 images, semivariogram value is plotted against lag distances. The direction of maximum correlation is determined to be $45^\circ$ with a ratio of 3 since the distance required to pass $\gamma_{\scriptsize \mbox{crit}}=0.9$ was largest in that direction and that distance is 3 times greater than the distance needed in the $-45^\circ$ direction. The directional semivariograms for the rotated and scaled field are shown on the right with a ratio of 1. }
\label{dir_variog}
\end{figure}

We can apply this same change of variables (\ref{anis_coord_transform}) to any two-dimensional spatial field that exhibits geometric anisotropy to achieve isotropy. For example, if we begin with a spatial field that exhibits its larger correlation length in the $45^\circ$ direction with $\tau=\ell_1/\ell_2=3$, the change of variables will rotate the field so the direction of maximum correlation length is in the $0^\circ$ direction and will then stretch the field along the new $y$-axis to remove the geometric anisotropy and create a new, isotropic field. This is shown in the middle in Figure \ref{dir_variog}. Once the spatial field has been adjusted in this way, a semivariogram can be fit to the transformed field as in the usual, isotropic case.

In order to adjust the spatial field to satisfy the isotropy assumptions in the way described above, we must ascertain $\theta$, the direction of maximum correlation length measured from the $x$-axis, and $\tau$, the ratio of the correlation length in the direction of $\theta$ to the correlation length in the direction orthogonal to $\theta$. Both of these parameters can be estimated using directional empirical semivariograms. Directional semivariograms are fit in a similar way as omnidirectional semivariograms in (\ref{empiricalSV}), but instead of taking all points separated by a distance $r$, we restrict the pairs of points to a certain angle, $\psi$. If we think of $\boldsymbol{w}_{i}$ and $\boldsymbol{w}_{j}$ as vectors, then $\psi$ is equivalent to the angle between $\boldsymbol{w}_{i}-\boldsymbol{w}_{j}$ and the $x$-axis. For example, if $\psi=0$, we restrict to all pairs of locations $w_i$ and $w_j$ on the same horizontal line, i.e., with the same $y$-coordinate. Formally, the empirical directional semivariogram can be defined as 
\begin{equation}
\label{directionalSV}
\hat{\gamma}_{\psi}(r)=\frac{1}{2n(r,\psi)}\sum_{(i,j)\left|\Vert\boldsymbol{w}_{i}-\boldsymbol{w}_j\Vert=r\right.,\phi_{ij}=\psi}\hspace{-1cm}[z(\boldsymbol{w}_i)-z(\boldsymbol{w}_j)]^2,
\end{equation}
where $\phi_{ij}$ denotes the angle between $\boldsymbol{w}_i-\boldsymbol{w}_j$ and the $x$-axis, and $n(r,\psi)$ is the number of points that are separated by a distance $r$ with angle of separation equal to $\psi$. It is common to calculate a directional semivariogram for $-90^\circ<\psi\le90^\circ$ in steps of either $15^\circ$ or $30^\circ$. We take a step size of $15^\circ$ here, which will result in 12 directional semivariograms.

Once the directional semivariograms have been calculated for each of the 12 different $\psi$ angles, we fit a common scatterplot smoother, the loess curve \cite{jacoby2000loess}, to the semivariogram values in each direction to achieve continuous curves. Then, to determine the ratio of correlation lengths, we can select a constant $\gamma_{\scriptsize\mbox{crit}}$ value between the nugget and sill and observe the distance required for the loess curve to surpass the height of $\gamma_{\scriptsize \mbox{crit}}$. The direction of maximum correlation, $\theta$, will require a larger distance to reach $\gamma_{\scriptsize\mbox{crit}}$ than other directions since the variance of the differences between values in that direction is expected to be smaller. The anisotropy ratio, $\tau$, can then be computed as the ratio between the distance in the direction of $\theta$ and the distance in the direction perpendicular to $\theta$.

This process is illustrated in Figure \ref{dir_variog}. The directional semivariograms are shown for the original, anisotropic field on the left. We can see that the correlation length is largest in the $45^\circ$ direction since the distance of $0.1684$ that it takes for the curve to pass $\gamma_{\scriptsize \mbox{crit}}=0.9$ is the largest of any direction. The range distance in the $-45^\circ$ direction is $0.0561$ and so the ratio of those ranges is $\tau=0.1684/0.0561=3$.

We can then rotate the field clockwise by $45^\circ$ and stretch it in the direction of the new $y$-axis by a factor of $\tau=3$ to achieve an isotropic field, as is done in the middle of Figure  \ref{dir_variog}. The directional semivariograms for the new field are shown on the right in Figure \ref{dir_variog}. It now takes a distance of $0.1684$ for the variogram values to pass $\gamma_{\scriptsize \mbox{crit}}$ for each $\psi$ angle, which means the ratio has been reduced to one, as it should be for an isotropic field. It is not always the case that we can reduce the ratio of these range values down to one, but we can reduce it enough for the field to be considered approximately isotropic.

\begin{algorithm}[b]
\caption{The Semivariogram Method for MAP Estimation with Anisotropic Whittle-Mat\'ern Prior.}
\label{alg2}
0. Set $\boldsymbol{x}_\alpha=\bb$.\\
1. Estimate $\theta$ and $\tau$ by computing directional semivariograms for $\boldsymbol{x}_\alpha$.\\
2. Transform the anisotropic spatial field coordinates, $\boldsymbol{w}$, to isotropic spatial field coordinates, $\boldsymbol{u}$, using (\ref{anis_coord_transform}).\\
3. Estimate $\thetab=(a_0,\sigma^2,\nu,\ell_1)$ by fitting an isotropic Mat\'ern semivariogram model to the transformed field. Then compute $\ell_2=\ell_1/\tau$.\\
4. Define the prior precision matrix, $\Pb$, by (\ref{equ:P}) using $\nu$, $\ell_1$, $\ell_2$, and $\theta$, compute $\alpha$ using (\ref{equ:GCV}), and compute $\boldsymbol{x}_\alpha$ using (\ref{equ:MAP2}).\\
5. Return to step 1 and repeat until $\theta$, $\tau$, $\nu$, $\ell_1$, and $\ell_2$ stabilize.
\end{algorithm}

Once we have obtained $\theta$ and $\tau$ and have changed the coordinates of the field, we can fit an isotropic omnidirectional semivariogram to estimate $\nu$ and $\ell_1$. Then we let $\ell_2=\ell_1/\tau$. All hyperparameters for use in (\ref{disc_anis}) will have been estimated and we can update these estimates iteratively using Algorithm \ref{alg2}. 
%\noindent{\bf Algorithm:} {\em The Semivariogram Method for MAP Estimation with Anisotropic Whittle-Mat\'ern Prior:}
%\begin{enumerate}
%\item[0.] Set $\boldsymbol{x}_\alpha=\bb$.
%
%\item[1.] Estimate $\theta$ and $\tau$ by computing directional semivariograms for $\boldsymbol{x}_\alpha$.
%
%\item[2.] Transform the anisotropic spatial field coordinates, $\boldsymbol{w}$, to isotropic spatial field coordinates, $\boldsymbol{u}$, using (\ref{anis_coord_transform}).
%
%\item[3.] Estimate $\thetab=(a_0,\sigma^2,\nu,\ell_1)$ by fitting an isotropic Mat\'ern semivariogram model to the transformed field. Then compute $\ell_2=\ell_1/\tau$.
%
%\item[4.] Define the prior precision matrix, $\Pb$, by (\ref{equ:P}) using $\nu$, $\ell_1$, $\ell_2$, and $\theta$, compute $\alpha$ using (\ref{equ:GCV}), and compute $\boldsymbol{x}_\alpha$ using (\ref{equ:MAP2}).
%
%\item[5.] Return to step 1 and repeat until $\theta$, $\tau$, $\nu$, $\ell_1$, and $\ell_2$ stabilize.
%\end{enumerate}
The convergence criteria for these hyperparameters are as follows: $\theta_j-\theta_{j-1}=0$, $\nu_j-\nu_{j-1}=0$, $|\ell^j_1-\ell_1^{j-1}|/\ell_1^{j-1}<0.01$, and $|\ell^j_2-\ell_2^{j-1}|/\ell_2^{j-1}<0.01$ where $\theta_j, \nu_j, \ell_1^j$ and $\ell_2^j$ denotes the $j$th iteration of the respective hyperparameter.

\subsection{Numerical Experiments}

\begin{figure}
\centering
\raisebox{2.3mm}{\includegraphics[width=1.78in]{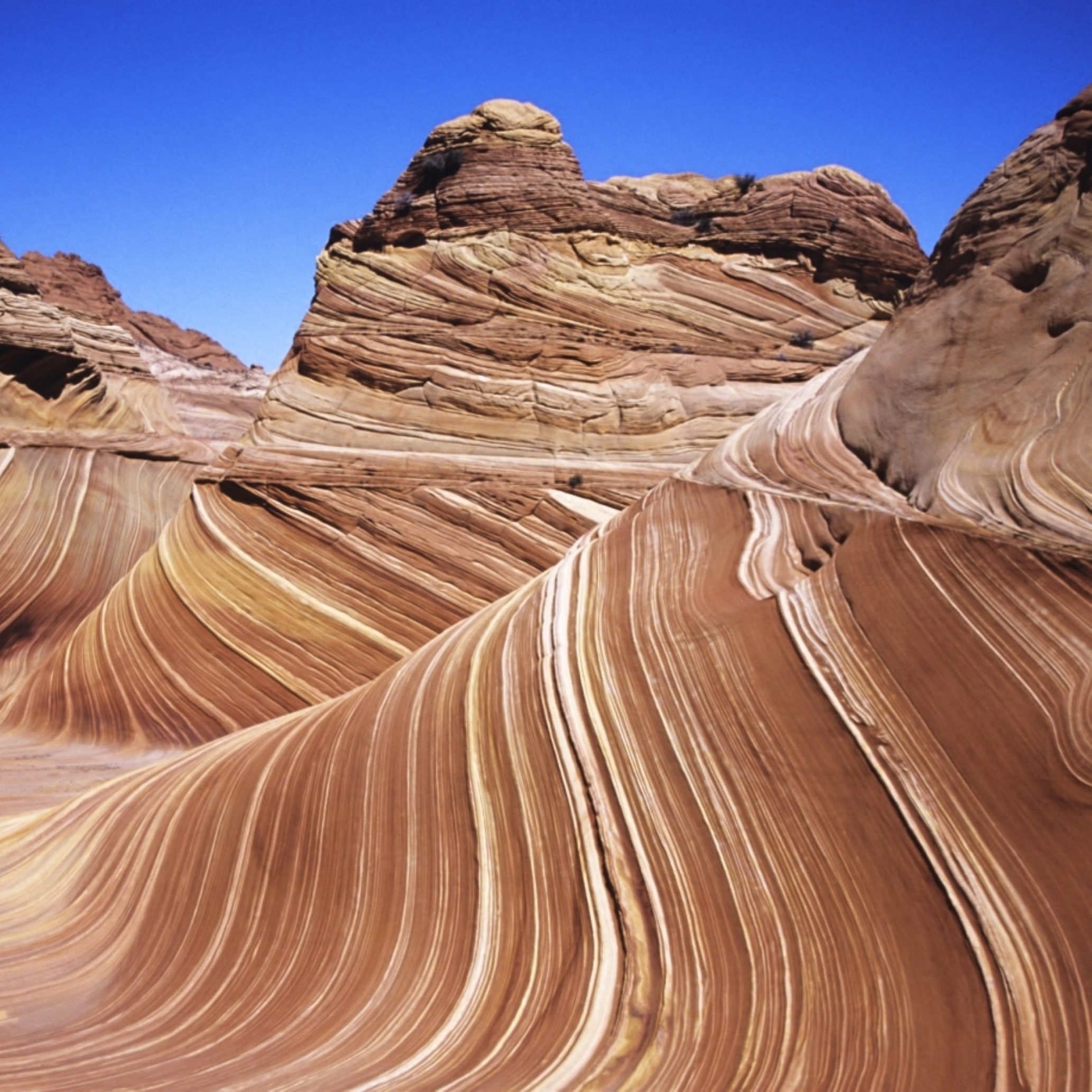}}\quad
\includegraphics[width=1.9in]{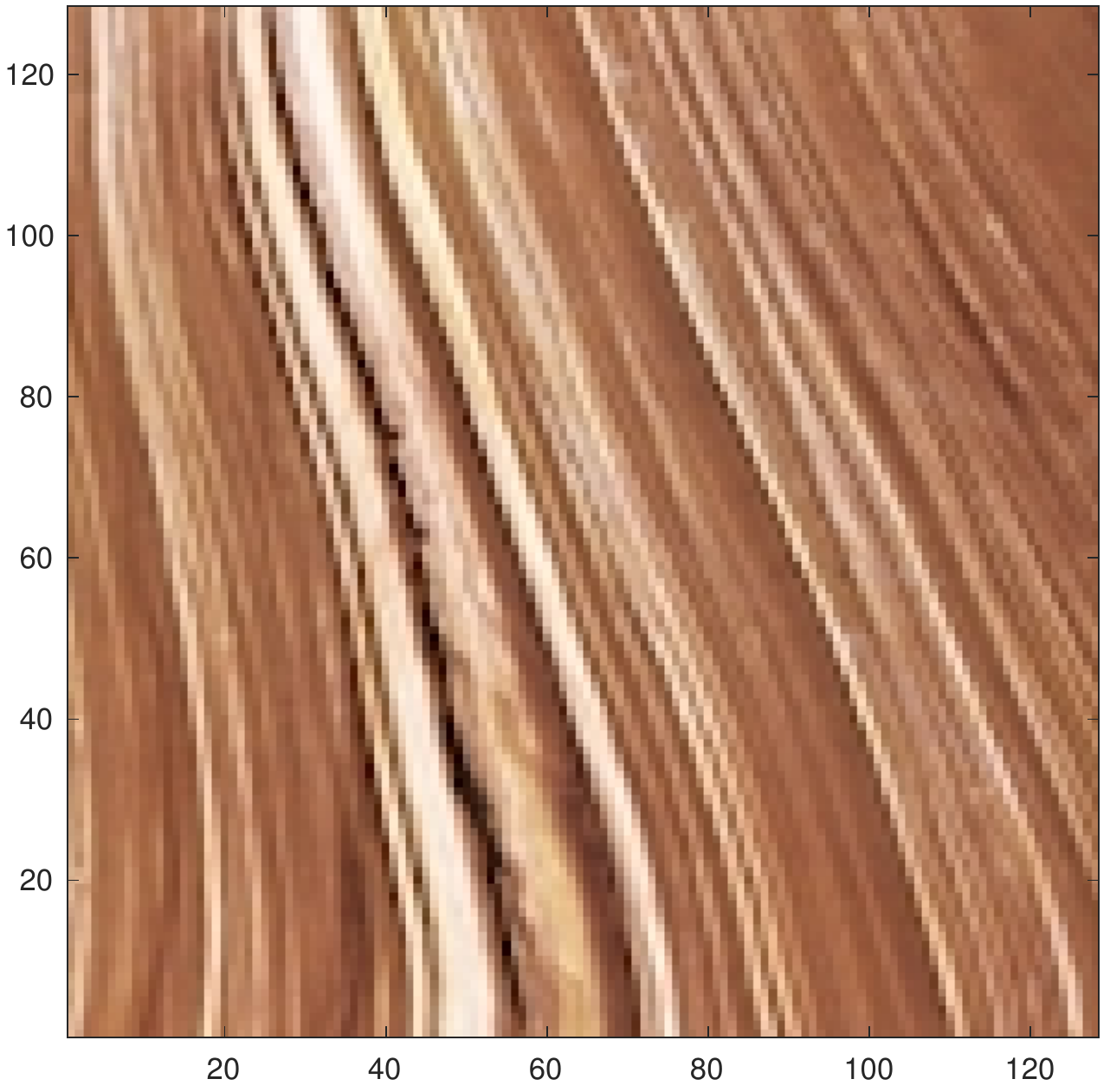}\quad
\includegraphics[width=1.9in]{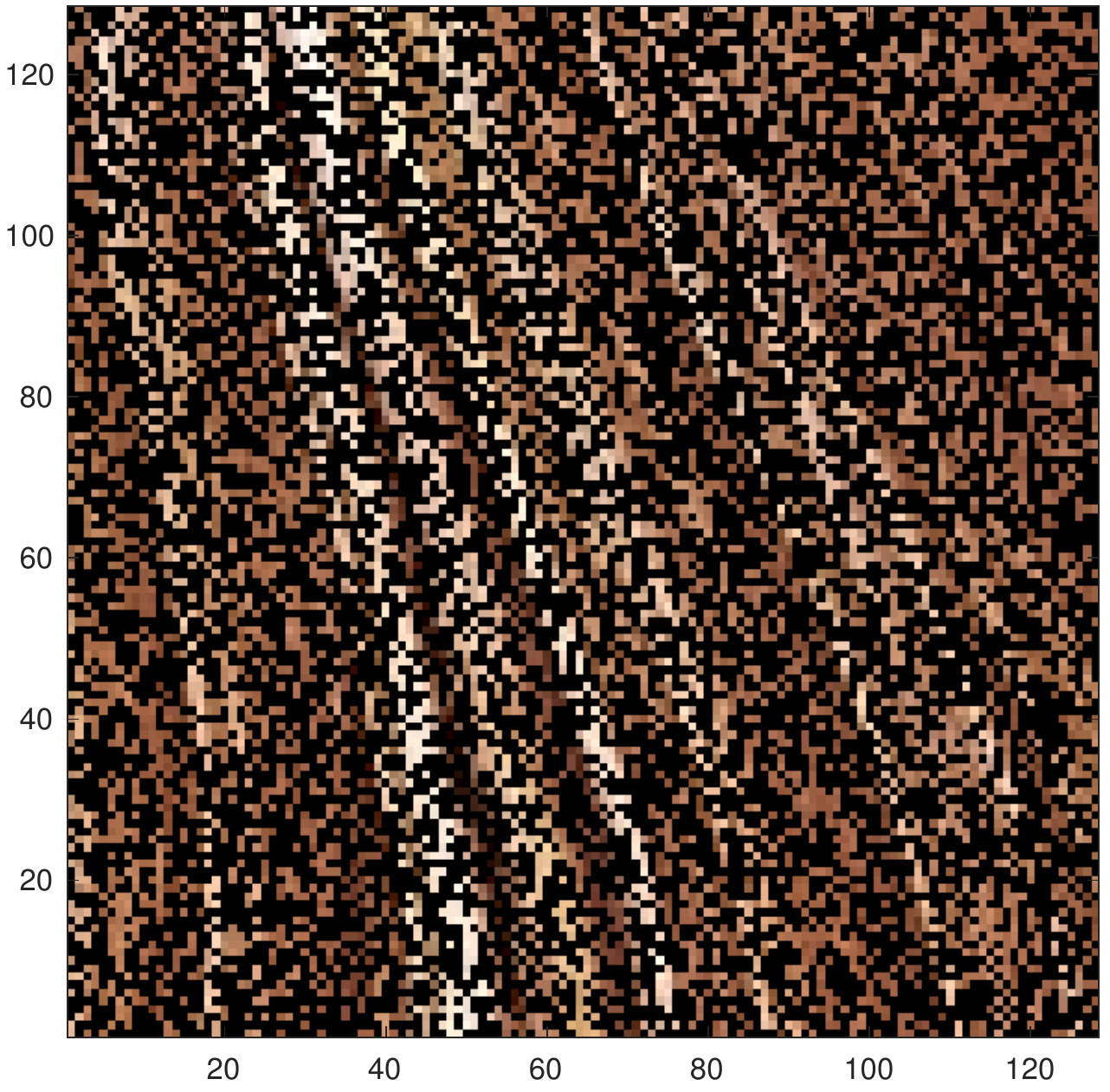}
\caption{Inpainting example. The original image showing the rock layers of the Wave in northern Arizona is given on the left. The true image used in the inpainting example is given in the middle. The masked image is given on the right.}
\label{inpaint_ex}
\end{figure}

We will illustrate the semivariogram method in the anisotropic case with a two-dimensional inpainting example. The original image, given on the left in Figure \ref{inpaint_ex}, shows a rock formation in Northern Arizona known as the Wave \cite{Arizona_Wave} where the layers of sandstone strata are clearly visible. We selected a subsection in the lower-middle of the image, shown in the middle of Figure \ref{inpaint_ex}, to illustrate our method. This will be the true image. We then added some noise and masked 60\% of the image. This is shown on the right in Figure \ref{inpaint_ex}. 

Like we saw in Section \ref{iso_results}, the prior will play a large role in the inpainting process since much of the image is missing. We will directly compare the solution using the anisotropic Whittle-Mat\'ern prior to the solution using the isotropic Whittle-Mat\'ern prior, both of which will have hyperparameters determined using semivariograms. Like before, the regularization parameter, $\alpha$, will be optimized using the highest correlation between the solution and the true image. 

After calculating the directional semivariograms for the image, the direction of maximum correlation was determined to be $-75^\circ$ for each color intensity. For the blue color-band, the correlation length in that direction was $\ell_1=0.1517$ and the correlation in the $15^\circ$ direction was $\ell_2=0.0101$, which gives a ratio of $\tau=15$. $\nu$ was determined to be $1$ and all of these hyperparameters converged in at most four iterations for each color and the initial $\theta$ estimate of $-75^\circ$ given in the first iteration remained unchanged throughout the process. When fitting an omnidirectional semivariogram to the masked image for the isotropic case, $\nu=2$ and $\ell=0.0142$.

\begin{figure}
\centering
\includegraphics[width=1.9in]{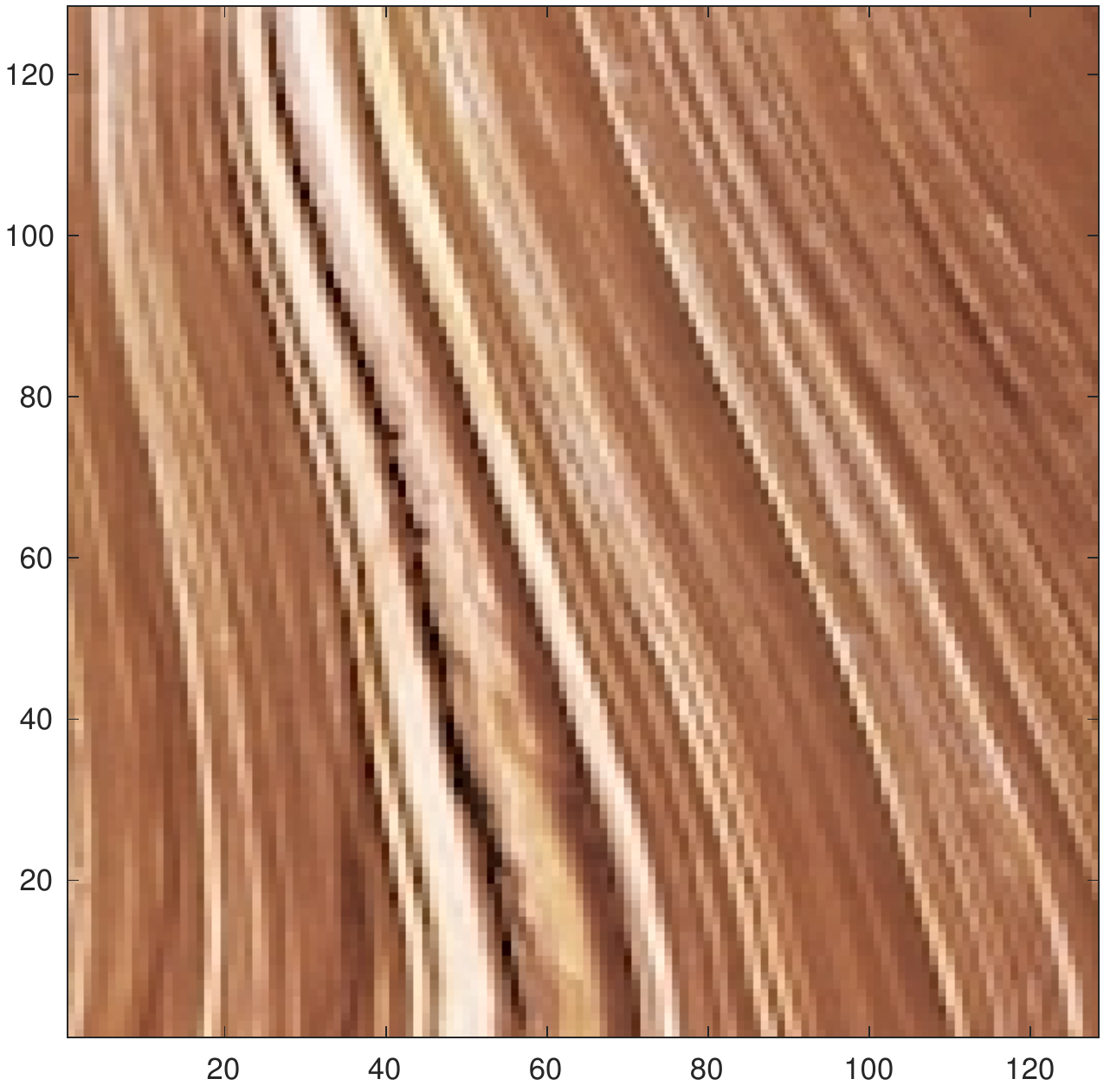}\quad
\includegraphics[width=1.9in]{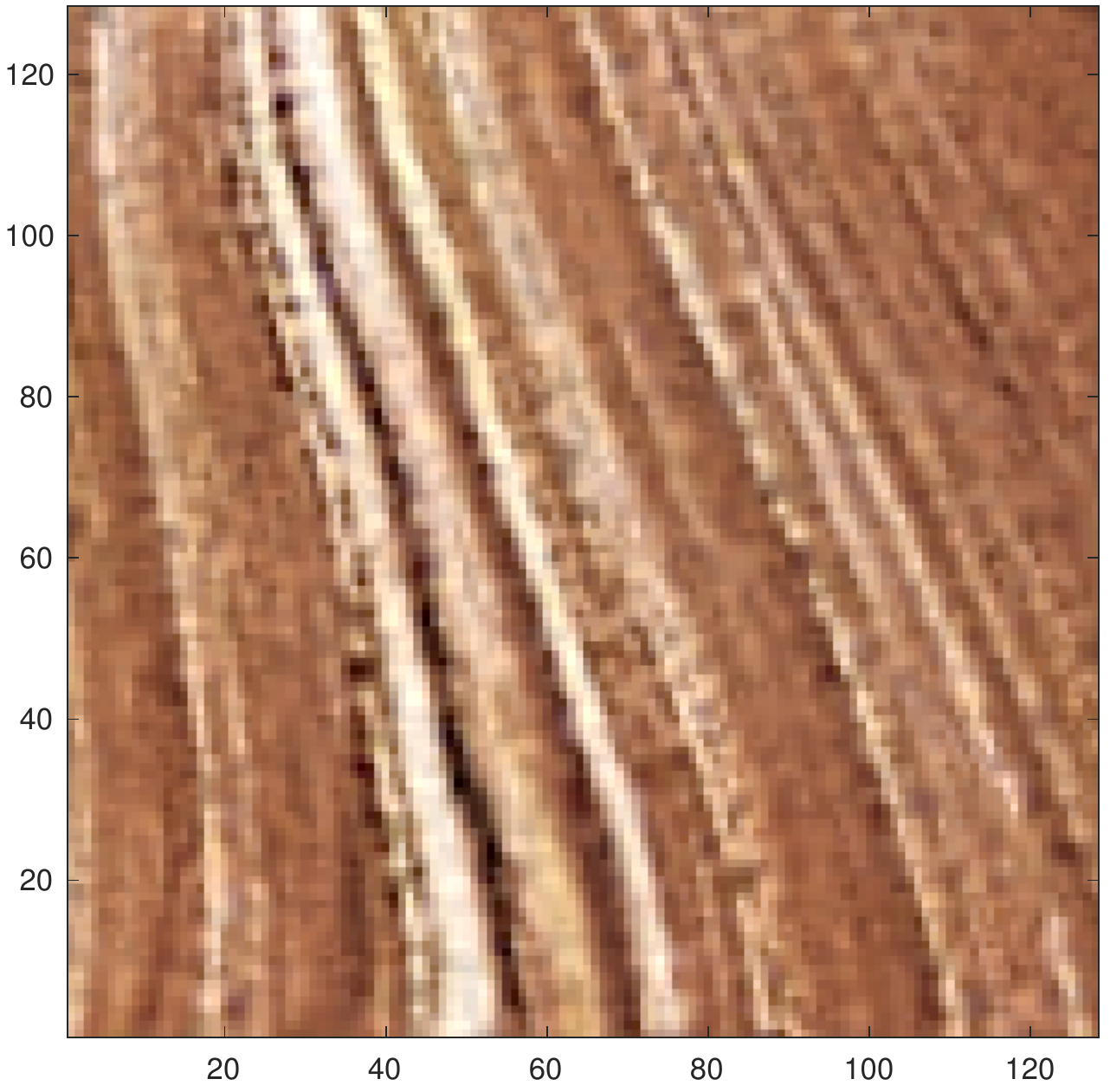}\quad
\includegraphics[width=1.9in]{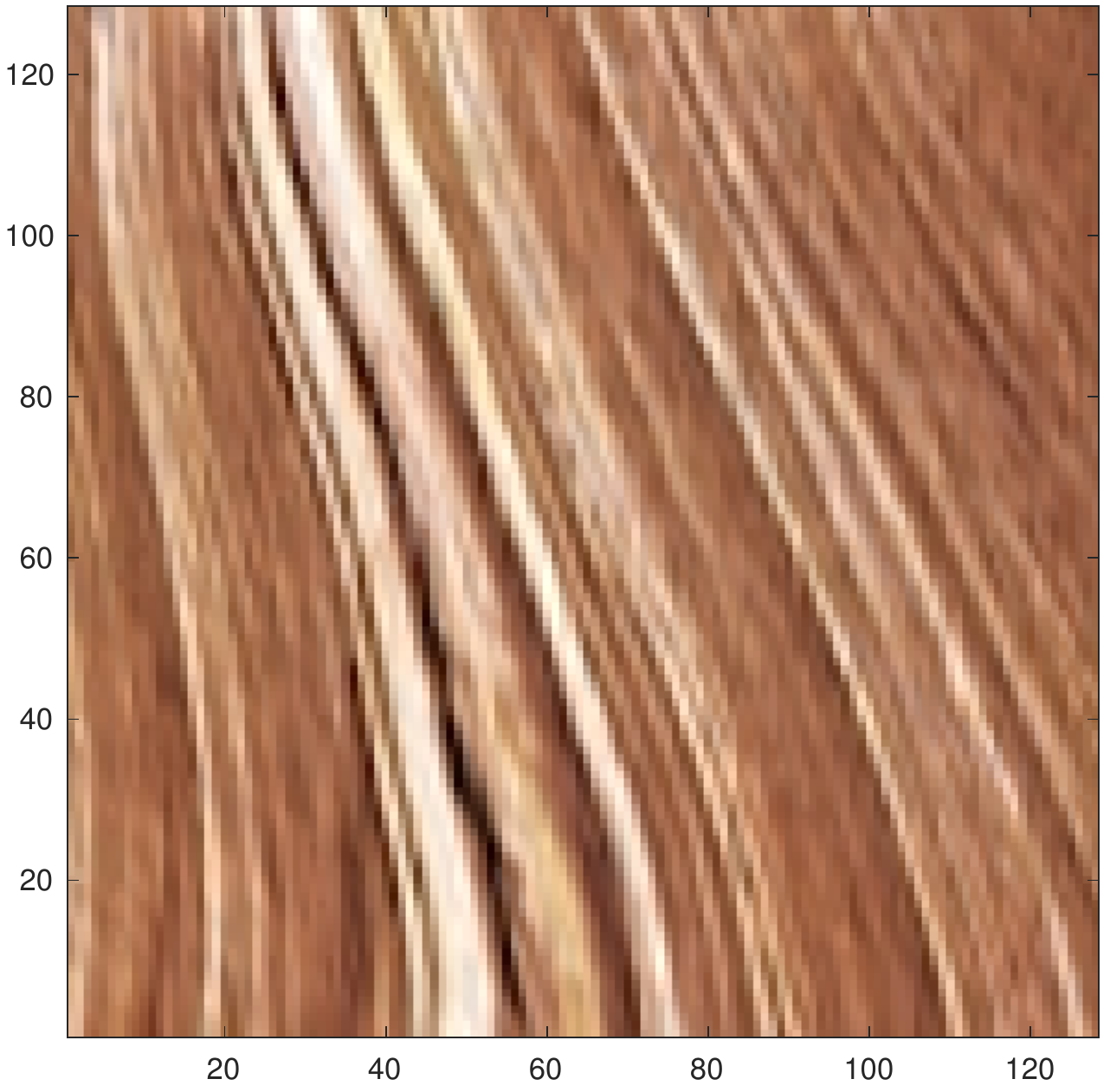}
\caption{Inpainting solutions. The true image (left) is given along with the the isotropic solution (middle) and the anisotropic solution (right).}
\label{inpaint_solutions}
\end{figure}

The reconstructions are given in Figure \ref{inpaint_solutions}. 
With the isotropic solution, the masking is removed, but since the prior assigns a very small correlation between each pixel, the reconstruction is noticeably spotty. The anisotropic solution, however, does a good job of removing the masking completely. The reconstruction is a bit smoother than the true image, but the original sandstone layers can be seen nicely. 

Some statistics of the reconstructions are given in Table \ref{table:inpaint_solutions}. Although the isotropic solution is still competitive, the anisotropic prior gives the reconstruction that most closely aligns with the true image. The isotropic solution has a mean absolute error (MAE) more than 57\% larger and a mean squared error (MSE) more than 180\% higher than those respective measures in the anisotropic case. The anisotropic reconstruction does fall short with the minimum value, however, which is farther from the truth than the solution given by the isotropic prior.
\begin{table}
\caption{\label{table:inpaint_solutions}Statistics for inpainting MAP estimates.} 
\begin{indented}
\lineup
\footnotesize
\item[]\begin{tabular}{@{}*{4}{l}}
\br                              
& True Image & Isotropic Covariance & Anisotropic Covariance \\
 \mr
$\bar{x}$ &  0.530 &  0.530 & 0.530  \\
s & 0.207  & 0.202 &  0.206 \\
Min & 0.000&  $\-0.053$& $\-0.065$\cr
$Q_1$ & 0.357 &  0.364 & 0.360 \cr
Median & 0.522  & 0.520 &  0.520 \cr
$Q_3$ & 0.678 &  0.676& 0.678\cr
Max & 1.000  & 1.112 & 1.093 \cr
$\rho_{\boldsymbol{x}_{\alpha},\boldsymbol{x}_{\scriptsize\mbox{true}}} $& &   0.944 & 0.981 \cr
%\mr
\mbox{Residual MAE} &  & 0.045 & 0.029  \cr
\mbox{Residual MSE} &  & 0.005 & 0.002 \cr
\br
\end{tabular}
\end{indented}
\end{table}

\subsection{Discussion}
Although the reconstruction with the anisotropic prior covariance matrix is better here, there are still some improvements that can be made. This example had a constant angle of maximum correlation length throughout the image and the ratio between maximum and minimum correlation was rather high, that is, greater than five. If either of these features fail to hold, the anisotropic prior often produces a reconstruction that performs slightly worse or offers no benefit over using an isotropic prior. We focus on the case when the angle of maximum anisotropy is not constant in the next section.

\section{Regional Anisotropy}
\label{sec:regional}
We have a way to define priors for isotropic and anisotropic spatial fields as long as that covariance structure is consistent for the entire field. In the case where the correlation length and angle of maximum anisotropy change throughout the image, we will want to model each of these regions with a different prior. 

\subsection{Regional Precision Matrix}

Suppose we have $k$ different regions in our image, each of which has a different covariance structure. We define $\Db_i$, $i=1,\dots,k$, as a masking matrix such that the only non-zero elements of $\Db_i\xb$ are those in region $i$. We will not allow for overlapping regions so that $ \sum_{i=1}^k \Db_k=\mathbf{I}$, the identity matrix. Now, to establish a prior for $\xb$ in this regional case, we take $\mbox{Cov}(\xb)=\mbox{Cov}(\Db_1\xb+\dots+\Db_k\xb)=\mbox{Cov}(\Db_1\xb)+\dots+\mbox{Cov}(\Db_k\xb)$, since each region is assumed independent due to not having any elements of $\xb$ in common. Define the best Whittle-Mat\'ern covariance structure, as chosen by a semivariogram, for region $i$ as $\Cb_i$ with corresponding precision matrix $\Pb_i=\Cb_i^{-1}$. Then $\mbox{Cov}(\Db_i\xb):=\Db_i\Cb_i\Db_i$. Thus, the prior for $\xb$ in this regional case has pdf
\begin{eqnarray}
\label{equ:regionalprior}
p(\boldsymbol{x}|\delta)\propto\exp\left(-\frac{\delta}{2}\boldsymbol{x}^T(\Db_1\Cb_1\Db_1+\dots+\Db_k\Cb_k\Db_k)^{-1}\boldsymbol{x}   \right),
\end{eqnarray}
which means our precision matrix is given by $\Pb=(\Db_1\Cb_1\Db_1+\dots+\Db_k\Cb_k\Db_k)^{-1}$. Note that (\ref{equ:regionalprior}) reduces to (\ref{equ:generalprior}) with $\Pb=\Pb_1$ when $k=1$. In general, the $\Cb_i$ matrices and $\Pb$ are dense, so actually constructing this precision matrix is infeasible for large problems. Additionally, FFTs cannot be used since $\Db_1\Cb_1\Db_1+\dots+\Db_k\Cb_k\Db_k$ is not circulant even if each $\Cb_i$ is. Thus, we seek an alternative expression such that the matrix-vector multiplication $\Pb\xb$ is achievable.

Without loss of generality, let $k=2$. Let 
\begin{eqnarray*}
\Cb_1 = \left[\begin{array}{cc} C_{1A} & C_{1B} \\C_{1C} & C_{1D} \end{array}\right]=\Pb_1^{-1} =\left[\begin{array}{cc} P_{1A} & P_{1B} \\P_{1C} & P_{1D} \end{array}\right]^{-1}
\end{eqnarray*}
and
\begin{eqnarray*}
\Cb_2= \left[\begin{array}{cc} C_{2A} & C_{2B} \\C_{2C} & C_{2D} \end{array}\right]=\Pb_2^{-1} =\left[\begin{array}{cc} P_{2A} & P_{2B} \\P_{2C} & P_{2D} \end{array}\right]^{-1}.
\end{eqnarray*}
Also assume that the regions are defined in a way that divides the region vertically (an assumption we will drop later) so that
\begin{eqnarray*}
\fl \Cb=\mbox{Cov}(\xb) = \mbox{Cov}(\Db_1\xb+\Db_2\xb) = \mbox{Cov}(\Db_1\xb)+\mbox{Cov}(\Db_2\xb)= \left[\begin{array}{cc} C_{1A} & 0 \\0 & C_{2D}\end{array}\right],
\end{eqnarray*}
which means our precision matrix is 
\begin{eqnarray*}
\Pb=\Cb^{-1} =\left[\begin{array}{cc} C_{1A}^{-1} &0 \\0 & C_{2D}^{-1} \end{array}\right].
\end{eqnarray*}

Using the block matrix inversion identity, it can be shown that $C_{1A}^{-1}= P_{1A}-P_{1B}P_{1D}^{-1}P_{1C}$ and $C_{2D}^{-1} =  P_{2D}-P_{2C}P_{2A}^{-1}P_{2B}$ and thus
\begin{eqnarray*}
\fl \Pb=\Cb^{-1} &= \left[\begin{array}{cc} C_{1A}^{-1} & 0\\0 & C_{2D}^{-1} \end{array}\right]=\left[\begin{array}{cc} P_{1A}-P_{1B}P_{1D}^{-1} P_{1C}& 0\\ 0 & P_{2D}-P_{2C}P_{2A}^{-1}P_{2B} \end{array}\right],
\end{eqnarray*}
which can be equivalently written as 
\begin{eqnarray*}
\fl \Pb&=\Db_1\Pb_1\Db_1 - \Db_1\Pb_1(\Db_2\Pb_1\Db_2)^\dagger \Pb_1\Db_1 + \Db_2\Pb_2\Db_2 - \Db_2\Pb_2(\Db_1\Pb_2\Db_1)^\dagger \Pb_2\Db_2.
\end{eqnarray*}
In general, for $k>2$, 
\begin{eqnarray}
\label{regional_precision}
\Pb&=\Cb^{-1} =(\Db_1\Cb_1\Db_1+\Db_2\Cb_2\Db_2+\dots+\Db_k\Cb_k\Db_k)^{-1}\nonumber\\
& = \sum_{i=1}^k\left(\Db_i\Pb_i\Db_i - \Db_i\Pb_i\Big[(\mathbf{I}_N-\Db_i)\Pb_i(\mathbf{I}_N-\Db_i)\Big]^\dagger \Pb_i\Db_i\right),
\end{eqnarray}
which, since each $\Pb_i$ is sparse, involves only sparse matrices. It is straightforward to show (\ref{regional_precision}) holds even in the case where the regions do not divide the region vertically by performing a reordering of the indices of $\xb$.

Since we have an expression for $\Pb$, we can now discuss how to perform the multiplication $\Pb\xb$. This will be needed to perform an iterative inverse method such as conjugate gradient to obtain the MAP estimator. Since each $\Db_i$ and $\Pb_i$ is sparse, each matrix vector multiplication in (\ref{regional_precision}) is efficient except the ones involving pseudoinverses. We can, however, take advantage of the lower-rank structure of $\big[(\Ib_N-\Db_i)\Pb_i(\Ib_N-\Db_i)\big]^\dagger$, which has rank $N-r_i$ where $r_i$ the rank of $\Db_i$. Let $\Pb_{i,\scriptsize{\mbox{nz}}}$ be the square matrix that consists of all rows and columns of $(\Ib_N-\Db_i)\Pb_i(\Ib_N-\Db_i)$ that have any nonzero elements. That is, keep row and column $j$ of $(\Ib_N-\Db_i)\Pb_i(\Ib_N-\Db_i)$ if $[\Ib_N-\Db_i]_{j,j}=1$. Then let $\Rb_i=\mbox{chol}(\Pb_{i,\scriptsize\mbox{nz}})$ such that $\Pb_{i,\scriptsize\mbox{nz}}=\Rb_i^T\Rb_i$ where chol denotes the Cholesky factorization and $\Rb_i$ is upper triangular. The Cholesky decomposition is known to be efficient for sparse, symmetric, positive definite matrices such as $\Pb_{i,\scriptsize\mbox{nz}}$ \cite{watkins2004fundamentals}. Then we can perform the multiplication of $\Db_i\Pb_i[(\Ib_N-\Db_i)\Pb_i(\Ib_N-\Db_i)]^\dagger \Pb_i\Db_i\xb$ in the following way:
\begin{enumerate}
\item[1.] Multiply $\boldsymbol{y}_i=\Pb_i(\Db_i\xb)$.
\item[2.] Extract the $N-r_i$ elements of $\boldsymbol{y}_i$ that correspond to the nonzero diagonal elements of $\Ib_N-\Db_i$: $\boldsymbol{y}_i(ind)$.
\item[3.] Define a variable $\boldsymbol{z}_i$ as an $N\times 1$ vector of zeros.
\item[4.] Multiply by $\Big((\Ib_N-\Db_i)\Pb_i(\Ib_N-\Db_i)\Big)^\dagger$ by taking $\boldsymbol{z}_i(ind)=\Rb_i\backslash(\Rb_i^T\backslash \boldsymbol{y}_i(ind))$.
\item[5.] Complete the multiplication $\Db_i(\Pb_i\boldsymbol{z}_i)$.
\item[6.] Repeat for $1\le i\le k$, so $\displaystyle \Pb\xb=\sum_{i=1}^k \Db_i(\Pb_i\boldsymbol{z}_i)$.
\end{enumerate}

Step 4 is the most costly since it requires both a forward and a backward substitution. This can be performed more efficiently for large regions since the rank of $(\Ib_N-\Db_i)\Pb_i(\Ib_N-\Db_i)$ is inversely related to the size of region $i$. Also, sparse reorderings, such as the symmetric approximate minimum degree permutation, can be used so $\Rb_i$ has fewer nonzero entries. 
The multiplication of $\Pb\xb$ must be performed for each iteration of CG, but each $\Rb_i$ can be stored ahead of time so the Cholesky decompositions need only be performed once. We saw some improvements in the performance of the CG algorithm when a preconditioner was used. The total number of iterations was approximately 21\% lower, which corresponded to about a 15\% overall time saving.

\subsection{Numerical Experiments}
\begin{figure}
\centering
\includegraphics[width=1.9in]{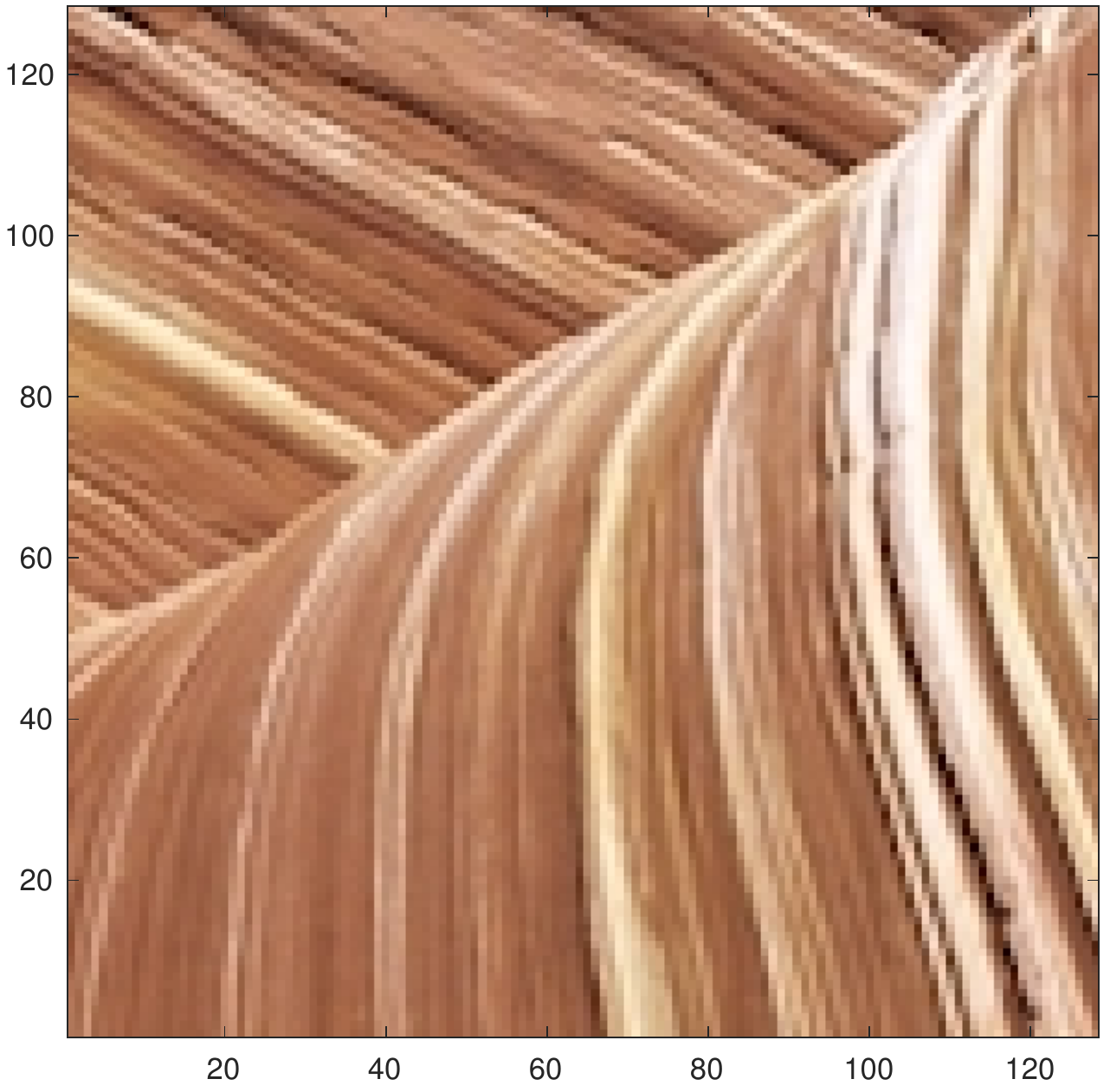}\quad
\raisebox{2.3mm}{\includegraphics[width=1.78in]{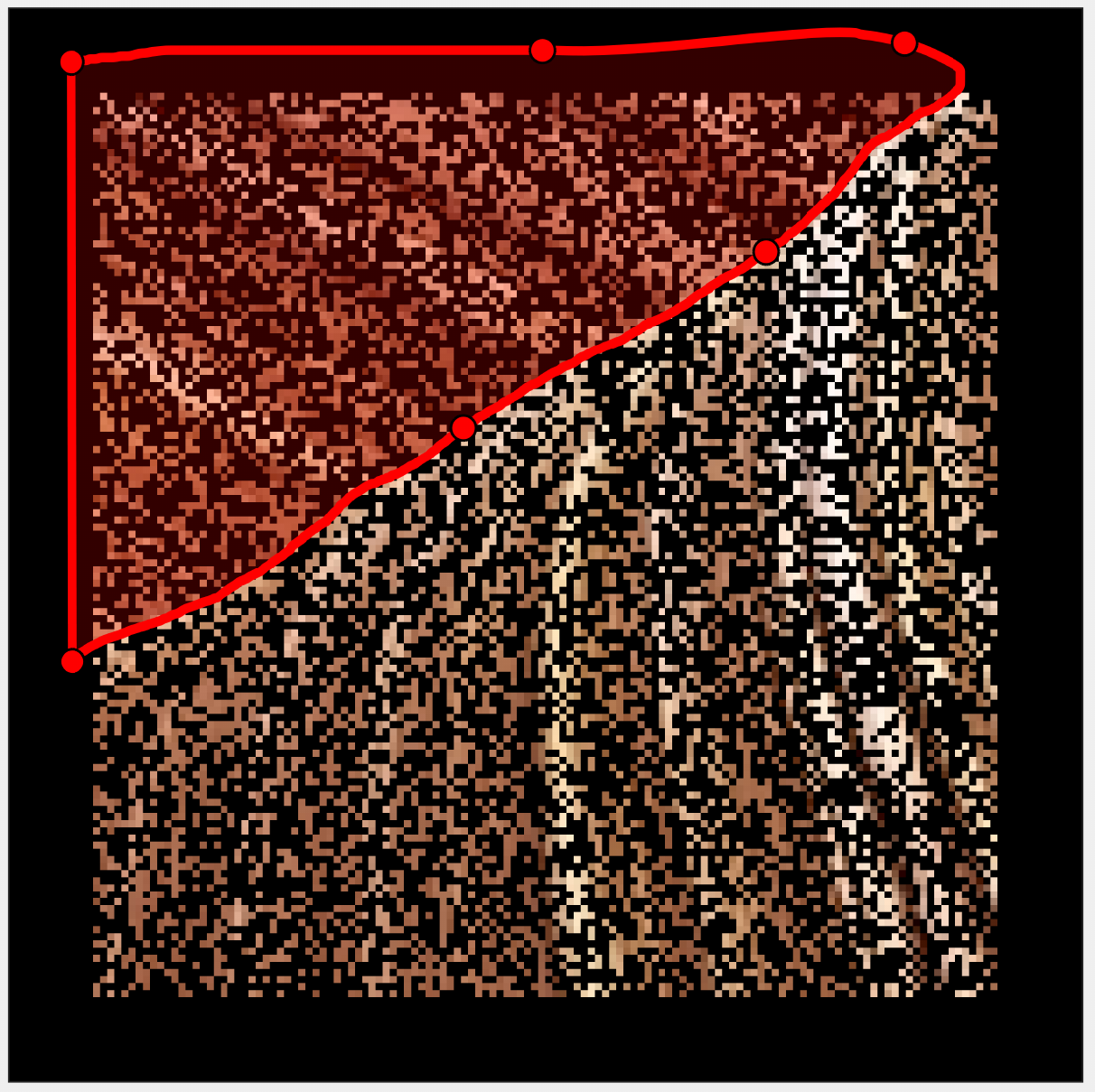}}\\
\includegraphics[width=1.9in]{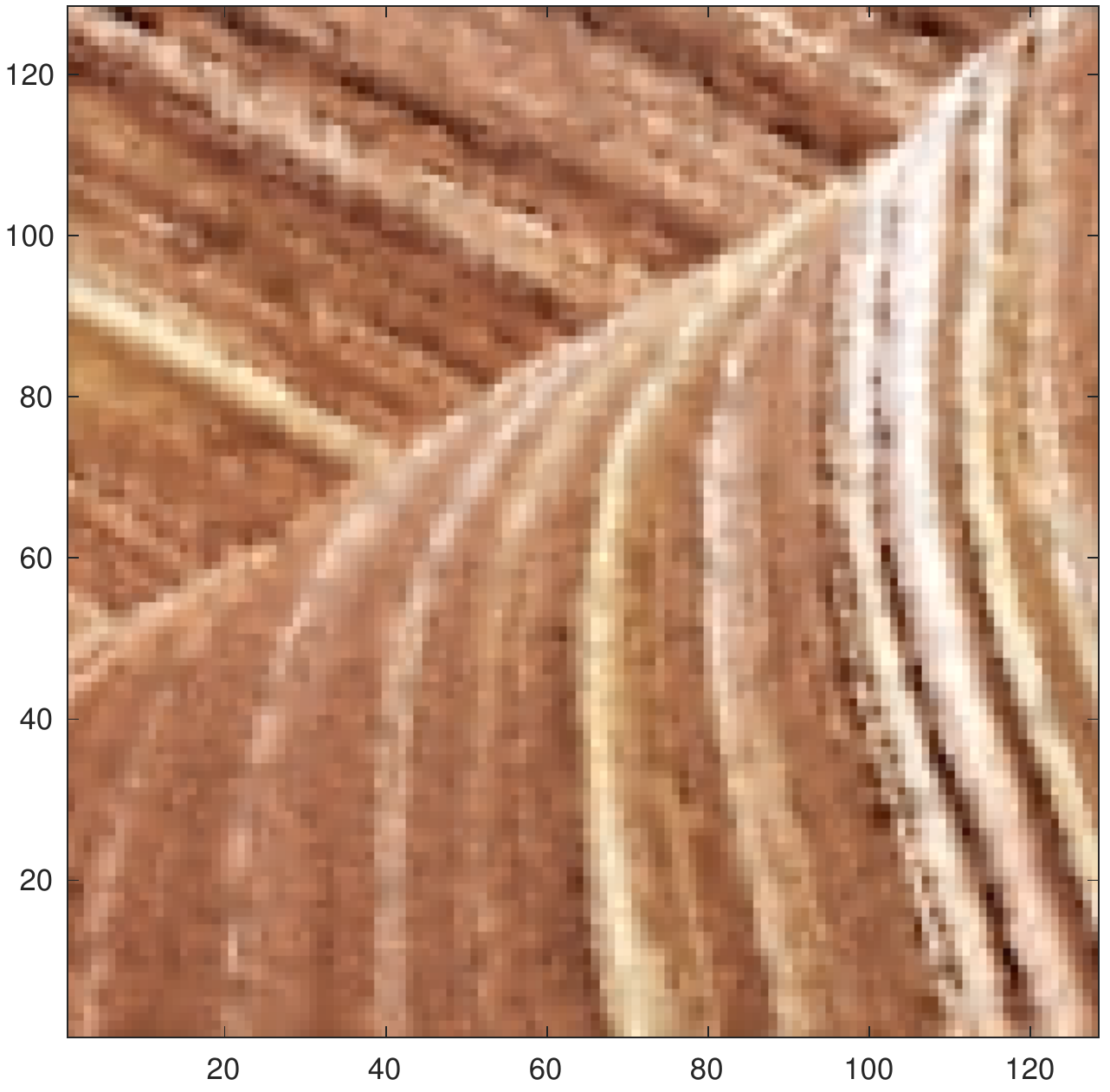}\quad
\includegraphics[width=1.9in]{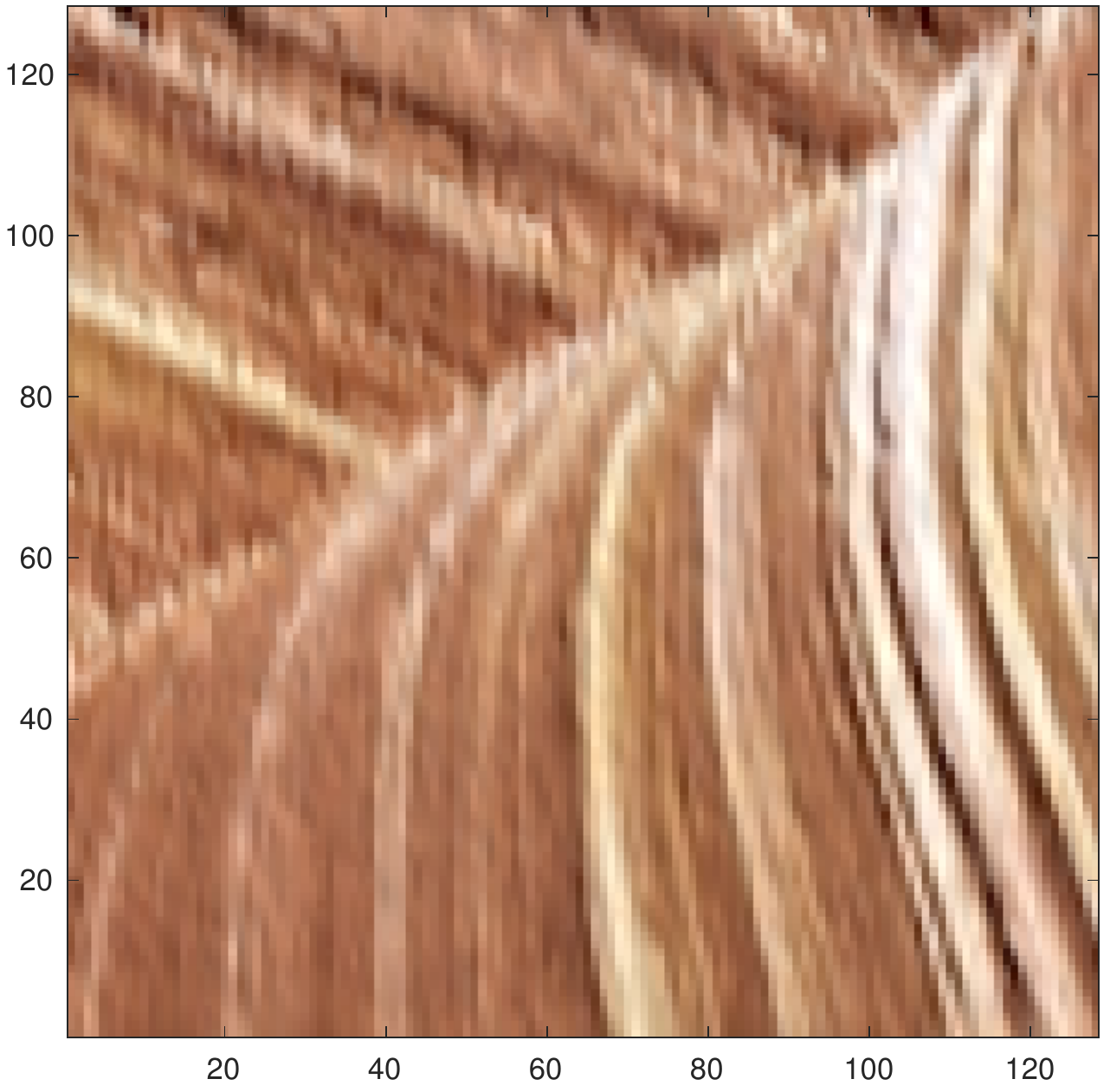}\quad
\includegraphics[width=1.9in]{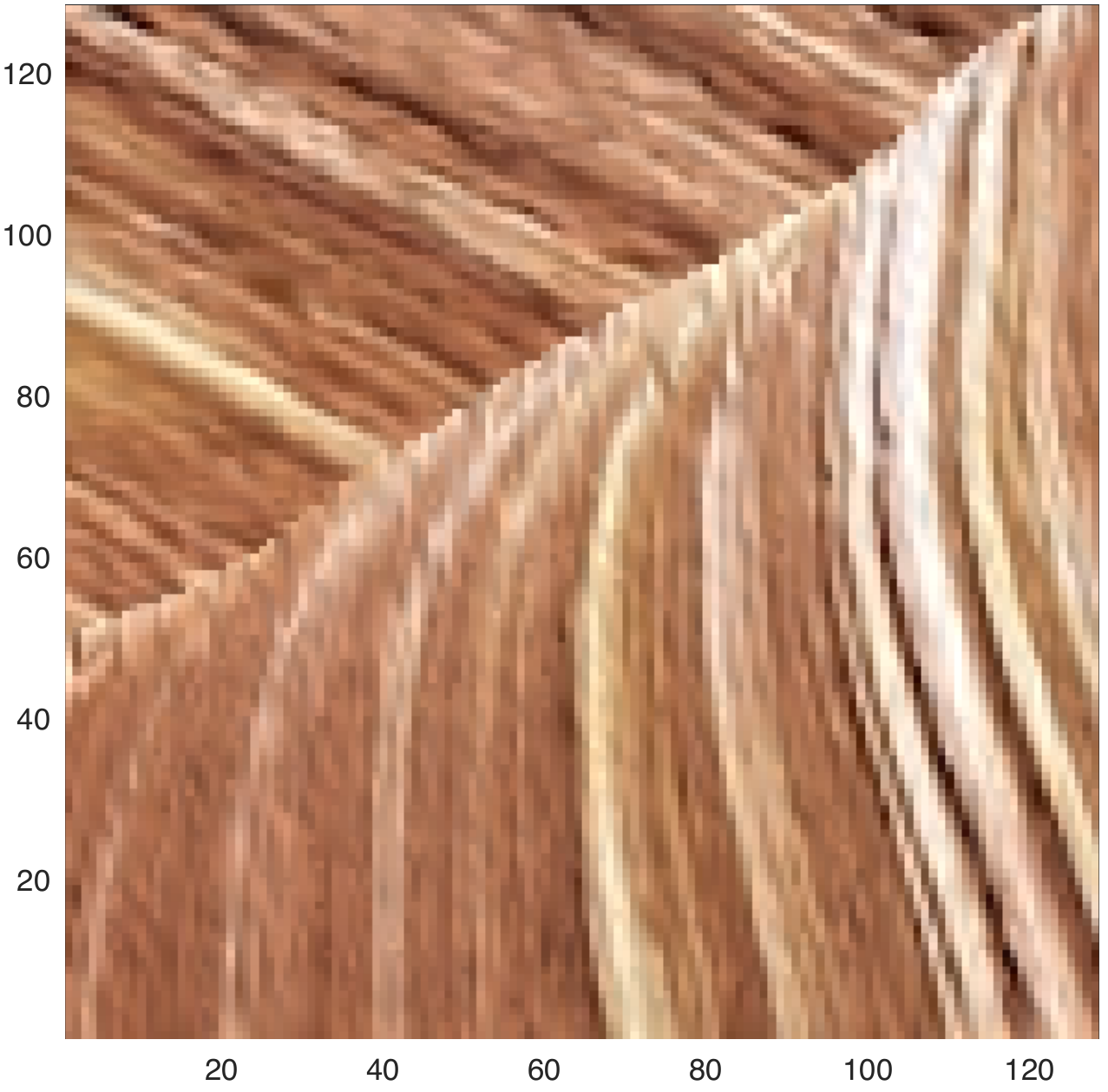}\quad
\caption{Inpainting solutions. The true image (top-left) is given along with the masked image (top-right), the isotropic solution (bottom-left) and the anisotropic solution (bottom-middle), and the regional solution (bottom-right).}
\label{inpaint_solutions2}
\end{figure}

We now consider an example where the angle of maximum anisotropy changes throughout the image. We take the central portion of the Wave image from Figure \ref{inpaint_ex} and again mask it so that 60\% of the image is blank. Then we attempt to inpaint the image using an isotropic prior, an anisotropic prior, and a regional anisotropic prior. The results are shown in Figure \ref{inpaint_solutions2}. The top-right image shows the masked picture as well as how the regions were chosen. The first region is shown with the red overlay while the second region is the remainder of the image. Semivariograms were fit to both regions and the top region was given a prior with an angle of maximum anisotropy of $-30^\circ$ while $\theta=-75^\circ$ for the bottom region. In the anisotropic solution given in the bottom-middle of the figure,  $\theta=-75^\circ$ throughout the image. Qualitatively, the regional solution in the bottom-right of the figure looks best.

Turning to Table \ref{table:inpaint_solutions2}, we can see the statistics comparing the different reconstructions. The isotropic and anisotropic solutions were similar in terms of the correlation and mean errors, but the regional solution is better in both of those categories and is similar in the others.

\begin{table}
\caption{\label{table:inpaint_solutions2}Statistics for regional inpainting MAP estimates.} 
%\begin{indented}
\lineup
\footnotesize
\begin{tabular}{@{}*{5}{l}}
\br                              
& True Image & Isotropic Covariance & Anisotropic Covariance & Regional Covariance \\
\mr
$\bar{x}$ & 0.567  & 0.565 &   0.566 &  0.566\cr
s & 0.207 &  0.200 &  0.202 & 0.207\\
Min & 0.000 & $\-0.014$ & $\-0.082$ & $\-0.032$  \cr
$Q_1$ & 0.400 & 0.402 & 0.402 & 0.398 \cr
Median & 0.565  & 0.564 & 0.564 & 0.562\cr
$Q_3$ & 0.722 & 0.717 & 0.718 & 0.720\cr
Max & 1.000   & 1.085 & 1.077 & 1.160 \cr
$\rho_{\boldsymbol{x}_{\alpha},\boldsymbol{x}_{\scriptsize\mbox{true}}} $& & 0.954 & 0.954 & 0.969\cr
%\mr
\mbox{Residual MAE} &  & 0.042 & 0.041 & 0.035 \cr
\mbox{Residual MSE} & & 0.004 & 0.004 &  0.003\cr
\br
\end{tabular}
%\end{indented}
\end{table}

\subsection{Discussion}
The regional covariance solution performed better in this example, but it does have some shortcomings. Firstly, it is best used when the distinction between regions is high. This is because the transition between regions when using this prior is abrupt, rather than smooth. Smoothing the transition between regions is something we leave to future work. Additionally, since multiplying $\Pb$ by $\xb$ requires inverting a matrix, this method can be slow when that matrix is large, which corresponds to a small region. Therefore, we suggest using small regions only when necessary. Alternatively, it is possible to solve a different inverse problem for each region independently and then combine the results. This will allow FFTs to be used since the precision matrix for each inverse problem will be in the form of (\ref{equ:P}).

\section{Conclusion}
\label{sec:conclusions}

In this paper, we introduced a method for selecting hyperparameters for use in the prior distribution of $\xb$ based on semivariogram modeling. We think of the noisy data as a spatial field and fit semivariograms to the noisy data and then iteratively to the MAP estimates to obtain point estimates for the prior hyperparameters. This method relies on the fact that the solution of the SPDE (\ref{equ:SPDE}) is a Gaussian process with zero mean and Mat\'ern covariance operator, which we have shown in detal. However, this connection requires an infinite domain, for us $\mathbb{R}^2$. For a finite domain, which is typically required for computations, the connection is broken, i.e., the SPDE solution is a zero mean Gaussian process without a Mat\'ern covariance operator. Fortunately, the connection can be restored by extending the finite computational domain. We showed how to systematically choose the extended domain using the Mat\'ern parameters.
The semivariogram method has the benefits of giving point estimates with a more intuitive interpretation while providing an objective way to choose an extension of the computational domain that is adequate for restoring the SPDE/Mat\'ern connection. We then applied the semivariogram method to an isotropic inpainting and deblurring example in two dimensions.

We generalized the isotropic results to the anisotropic case and showed the semivariogram method can be applied as well by using directional semivariograms and the anisotropic SPDE (\ref{SPDE_anis}). An inpainting example comparing reconstructions using isotropic and anisotropic priors was presented. Finally, we discussed an even more general case when the image has regions with differing correlation lengths and angles of maximum correlation, which requires a sparse precision matrix that can be obtained via a discretized SPDE for each region. One more example was shown that yielded good solutions.

\ack
J. Bardsley acknowledges support from the Gordon Preston Fellowship offered by the School of Mathematics at Monash University. T. Cui acknowledges support from the Australian Research Council, under grant number CE140100049 (ACEMS). We would also like to acknowledge the assistance of Dr. Jon Graham at the University of Montana with the semivariogram methodology.

\section*{References}
\bibliographystyle{unsrt}  
\bibliography{manuscript_IP_final}

\end{document}